\documentclass{gtart_h}

\def\ifplaintex{\expandafter\ifx\csname documentclass\endcsname\relax}

\def\gtp{{\mathsurround=0pt\it $\cal G\mskip-2mu$eometry \&\ 
$\cal T\!\!$opology $\cal P\!$ublications}}  

\def\recd{{\small Received:\qua\receiveddate\ifx\reviseddate\relax
\else\qquad Revised:\qua\reviseddate\fi\par}} 

\def\lognumber#1{\def\thelognumber{#1}}
\def\volumenumber#1{\def\thevolumenumber{#1}}
\def\volumeyear#1{\def\thevolumeyear{#1}}
\def\papernumber#1{\def\thepapernumber{#1}}
\def\pagenumbers#1#2{\def\startpage{#1}\def\finishpage{#2}}
\def\published#1{\def\publishdate{#1}}

\def\received#1{\def\receiveddate{#1}}
\def\revised#1{\def\reviseddate{#1}}
\def\accepted#1{\def\accepteddate{#1}}

\def\asciiaddress#1{\def\theasciiaddress{#1}}
\def\asciiemail#1{\def\theasciiemail{#1}}

\long\def\asciiabstract#1{\long\def\theasciiabstract{#1}}

\let\\\par\let\thelognumber\relax\let\thevolumenumber\relax
\let\thepapernumber\relax\let\thevolumeyear\relax\let\startpage\relax
\let\finishpage\relax\let\publishdate\relax\let\receiveddate\relax
\let\reviseddate\relax\let\accepteddate\relax\let\theasciititle\relax
\let\theasciiauthors\relax\let\theasciiaddress\relax
\let\theasciiabstract\relax

\let\theasciiemail\relax

\ifplaintex
\font\logobig=cmssbx10 scaled 3836
\font\logomed=cmssbx10 scaled 2557
\else
\font\logobig=cmssbx10 scaled 4200
\font\logomed=cmssbx10 scaled 2800
\fi

\long\def\makeagttitle{   
\count0=\startpage
\agt\hfill      
\hbox to 45truept{\vbox to 0pt{\vglue -13truept{\logomed A\kern -.37em{\logobig 
T}\kern -.38em G}\vss}\hss}
\break
{\small Volume \thevolumenumber\ (\thevolumeyear)
\startpage--\finishpage\nl
Published: \publishdate}

\vglue .25truein

{\parskip=0pt\leftskip 0pt plus
1fil\def\\{\par\smallskip}{\Large\bf\thetitle}\par\medskip} \vglue
0.05truein

{\parskip=0pt\leftskip 0pt plus 1fil\def\\{\par}{\sc\theauthors}
\par\medskip}%
 
\vglue 0.03truein 

{\small\leftskip 25truept\rightskip 25truept{\bf Abstract}\stdspace\theabstract

{\bf AMS Classification}\stdspace\theprimaryclass
\ifx\thesecondaryclass\relax\else; \thesecondaryclass\fi\par
{\bf Keywords}\stdspace \thekeywords\par}\vglue 7truept

}   

\ifplaintex
\font\phead=cmsl9 scaled 950
\font\pnum=cmbx10 scaled 913
\font\pfoot=cmsl9 scaled 950
\headline{\vbox to 0pt{\vskip -4.5mm\line{\small\phead\ifnum
\count0=\startpage ISSN 1472-2739 (on-line) 1472-2747 (printed)
\hfill {\pnum\folio}\else\ifodd\count0\def\\{ }%
\ifx\theshorttitle\relax\thetitle\else\theshorttitle\fi\hfill{\pnum\folio}
\else\def\\{ and }{\pnum\folio}\hfill\ifx\theshortauthors\relax\theauthors
\else\theshortauthors\fi\fi\fi}\vss}}
\footline{\vbox to 0pt{\vglue 0mm\line{\small\pfoot\ifnum\count0=\startpage
\copyright\ \gtp\hfill\else
\agt, Volume \thevolumenumber\ (\thevolumeyear)\hfill\fi}\vss}}
\else
\font=cmsl9 scaled 1050
\font\lnum=cmbx10 
\font=cmsl9 scaled 1050
\makeatletter
\def\@oddhead{{\small\lhead\ifnum\count0=\startpage ISSN 1472-2739 
(on-line) 1472-2747 (printed)\hfill {\lnum\number\count0}\else\ifodd\count0
\def\\{ }\ifx\theshorttitle\relax \thetitle \else\theshorttitle\fi\hfill
{\lnum\number\count0}\else\def\\{ and }{\lnum\number\count0}
\hfill\ifx\theshortauthors\relax 
\theauthors\else\theshortauthors\fi\fi\fi}}\def\@evenhead{\@oddhead}
\def\@oddfoot{\small\lfoot\ifnum\count0=\startpage\copyright\ \gtp\hfill\else
\agt, Volume \thevolumenumber\ (\thevolumeyear)\hfill\fi}
\def\@evenfoot{\@oddfoot}
\makeatother
\fi
\let\maketitlepage\makeagttitle
\let\makeshorttitle\maketitlepage
\let\maketitle\maketitlepage

\newwrite\gtoutfile
\long\gdef\makeheadfile{  
{\def\\{, }\def\s{ }
\immediate\openout\gtoutfile head.xxx
\immediate\write\gtoutfile{Proxy-for: \ifx\theasciiauthors\relax
\theauthors\else\theasciiauthors\fi\s<\ifx\theasciiemail\relax\theemail\else\theasciiemail\fi>}
\immediate\write\gtoutfile{\noexpand\\}
\immediate\write\gtoutfile{Authors: \ifx\theasciiauthors\relax
\theauthors\else\theasciiauthors\fi}
{\def\\{ }\immediate\write\gtoutfile{Title: \ifx\theasciititle\relax
\thetitle\else\theasciititle\fi}}
\immediate\write\gtoutfile{Subj-class: GT or SG, GR etc}
\immediate\write\gtoutfile{MSC-class: \theprimaryclass\ifx\thesecondaryclass\relax\else, \thesecondaryclass\fi}
\immediate\write\gtoutfile{Journal-ref: Algebr. Geom. Topol. \thevolumenumber\s
(\thevolumeyear) \startpage-\finishpage}
\immediate\write\gtoutfile{Comments: Published by Algebraic and
Geometric Topology at}
\immediate\write\gtoutfile{\s\s\s  http://www.maths.warwick.ac.uk/agt/AGTVol\thevolumenumber/agt-\thevolumenumber-\thepapernumber.abs.html}
\immediate\write\gtoutfile{\noexpand\\}
\immediate\write\gtoutfile{}
\ifx\theasciiabstract\relax
\immediate\write\gtoutfile{\theabstract}\else
\immediate\write\gtoutfile{\theasciiabstract}\fi
\immediate\write\gtoutfile{}
\immediate\write\gtoutfile{\noexpand\\}
\immediate\write\gtoutfile{}
\immediate\closeout\gtoutfile}}  

\def\maketitlepage{\makeagttitle\makeheadfile}
\let\makeshorttitle\maketitlepage
\let\maketitle\maketitlepage

\lognumber{25}
\volumenumber{4}
\volumeyear{2004}
\papernumber{25}
\published{15 July 2004}
\pagenumbers{543}{569}
\received{26 February 2004}
\revised{8 July 2004}
\accepted{9 July 2004}

\usepackage{amsmath,amssymb,graphicx}
\usepackage{color}

\usepackage[all]{xy}

\newtheorem{thm}{Theorem}[section]

\newtheorem{lem}[thm]{Lemma}

\newtheorem{prop}[thm]{Proposition}
\theoremstyle{definition}
\newtheorem{defn}[thm]{Definition}

\newcommand{\Z}{\ensuremath{{\mathbb{Z}}}}
\newcommand{\R}{\ensuremath{{\mathbb{R}}}}

\newcommand{\SP}{\ensuremath{{\mathbb{S}}}}

\newcommand{\al}{\alpha}
\newcommand{\e}{\epsilon}
\newcommand{\s}{\sigma}
\newcommand{\Ga}{\ensuremath{{\Gamma}}}
\newcommand{\ga}{\ensuremath{{\gamma}}}
\newcommand{\Om}{\Omega}
\newcommand{\Si}{\Sigma}
\newcommand{\F}{\mathcal{F}}

\newcommand{\ul}{\underline{l}}

\newcommand{\rar}{\rightarrow}

\newcommand{\vs}{\vspace{0mm}} 

\newcommand{\nk}{\binom{n}{k}}
\newcommand{\del}{\partial}

\newcommand{\Aut}{\operatorname{Aut}}
\newcommand{\A}{\operatorname{A}}
\newcommand{\Htpy}{\operatorname{Htpy}}

\newcounter{samcounter}

\newcounter{samcounter2}
\newenvironment{sam2}[3]{
\begin{list}{#2#1{samcounter2}#3\stepcounter{samcounter2}}{\setcounter{samcounter2}{1}}}
{\end{list}}

\begin{document}

\title{Automorphisms of free groups with boundaries}

\authors{Craig A. Jensen\\Nathalie Wahl}
      \address{Department of Mathematics, University of 
New Orleans\\New Orleans, LA 70148, USA}
      \secondaddress{Department of Mathematics, University of Aarhus\\Ny 
Munkegade 116, 8000 Aarhus, DENMARK}
\asciiaddress{Department of Mathematics, University of New Orleans\\New 
Orleans, LA 70148, USA\\and\\Department of Mathematics, University of 
Aarhus\\Ny Munkegade 116, 8000 Aarhus, DENMARK}

\gtemail{\mailto{jensen@math.uno.edu}{\rm\qua and\qua}\mailto{wahl@imf.au.dk}}
\asciiemail{jensen@math.uno.edu, wahl@imf.au.dk}

\begin{abstract}
The automorphisms of free groups with boundaries form a family of groups $\A_{n,k}$ closely related to  mapping class groups,  
with the standard automorphisms of free groups as  $\A_{n,0}$ and (essentially) the symmetric automorphisms of free groups as $\A_{0,k}$. 
We construct a contractible space $L_{n,k}$ on which 
$\A_{n,k}$ acts with finite stabilizers  and finite quotient  space and 
deduce a range for the virtual cohomological dimension of $\A_{n,k}$. 
We also give a presentation of the groups and calculate their first homology group.
\end{abstract}
\asciiabstract{%
The automorphisms of free groups with boundaries form a family of
groups A_{n,k} closely related to mapping class groups, with the
standard automorphisms of free groups as A_{n,0} and (essentially) the
symmetric automorphisms of free groups as A_{0,k}.  We construct a
contractible space L_{n,k} on which A_{n,k} acts with finite
stabilizers and finite quotient space and deduce a range for the
virtual cohomological dimension of A_{n,k}.  We also give a
presentation of the groups and calculate their first homology group.}


\maketitle

\section{Introduction}

The group $A_{n,k}$ of automorphisms of free groups with $k$ boundaries was introduced in \cite{w} to 
show that the natural map from the stable mapping class group of surfaces to the stable automorphism group of free groups gives 
an infinite loop map on the classifying spaces of the groups after plus construction. 
The ``boundaries'' make it possible to define appropriate gluing operations.

There are several ways to define these groups.

Geometrically, $\A_{n,k}$ is the group of components of the homotopy equivalences of a certain genus $n+k$ graph $G_{n,k}$ (see Figure~\ref{Ggk})
fixing $k$ circles in the graph. The mapping class group $\Ga_{g,k+1}$ of a surface of genus $g$ 
with $k+1$ boundary components is in fact  the subgroup 
of $\A_{2g,k}$ fixing one more circle (see Section~\ref{Anksec}). 
An alternative geometric description is given in \cite{hw} in terms of the  mapping class groups of a 3-dimensional manifold.

Algebraically, we first define a related group $\A_n^k$, which corresponds to mapping class groups with punctures rather than boundaries. 
Let $F_{n+k}$ denote the free group on the $n+k$ generators $x_1,\dots,x_n,y_1,\dots,y_k$, and let $\Aut(F_{n+k})$ be its automorphism 
group. 
Then $\A_n^k:=\{\phi\in \Aut(F_{n+k})\,|\, \phi(y_i)\sim y_i\ \textrm{for}\ 1\le i\le k\}$, where $\sim$ means `conjugate to'. 
The groups $\A_n^k$ are thus  natural intermediates 
between $\A_n^0=\Aut(F_n)$ and $\A_0^k$ which is known as the pure symmetric automorphism group of $F_k$. 
The group $\A_{n,k}$ is a central extension
    $$1\to \Z^k \to \A_{n,k}\to \A_n^k \to 1$$
(see Section~2 for details).  

The symmetric automorphisms of free groups have been studied extensively, 
in particular for their relation to automorphism groups of free 
products (see for example \cite{bmmm,c,jmm,mc4,mm}) and for their relation to motion groups of strings 
\cite{bl,d,g,mm}; 
they give a generalization of the braid groups, which is the mapping class group of 
a genus 0 punctured surface. 
Our groups are the higher genus analogues, coming from surfaces with punctures for $\A_n^k$ or with boundaries for $\A_{n,k}$.

In \cite[Section~4]{l}, Levitt defines a class of groups, which he calls {\em relative automorphism groups},
 in order to study automorphism groups of hyperbolic groups and 
of graphs of groups. The groups $A_{n,k}$ are examples of relative automorphism groups.

\vs

In this paper, we construct a contractible space $L_{n,k}$ of dimension $2n+3k-2$
(provided $n \geq 1$; if $n=0$ it is of dimension $2k-1$)
on which the group $\A_{n,k}$ acts with finite stabilizers 
and such that the quotient space is a finite complex. The space $L_{n,k}$ is an analogue of the spine of Auter space $X_{n}$. 
Recall that $X_{n}$ is a space of graphs of genus $n$ marked by a rose $\vee_{n}S^1$. In $L_{n,k}$, we consider graphs of genus $n+k$ marked by the 
graph $G_{n,k}$ of Figure~\ref{Ggk}. 
We also construct an analogous space $L_n^k$ for the group $A_n^k$ so that there is a sequence of contractible spaces 
$$F_{n,k}\to L_{n,k} \to L_n^k$$ 
corresponding to the group extension $\Z^k\to A_{n,k}\to A_n^k$. The fiber $F_{n,k}$ is homeomorphic 
to $\R^k$ with the standard $\Z^k$-action.

We use $L_{n,k}$  to give the following range for the virtual cohomological dimension of $\A_{n,k}$: 
 $2n+2k-1\le vcd(\A_{n,k})\le 2n+3k-2 $ when $k \geq 1$
(Theorem~\ref{vcd}).

By the work of McCool \cite{mc3}, we know that $\A_n^k$ is finitely presented. We use his work to produce a presentation for $\A_n^k$ 
and $\A_{n,k}$ (Theorems~\ref{presAnk1} and \ref{presAnk2}). We recover the presentation he gave for $\A_0^k$ in \cite{mc4}.
Theorem~\ref{H1} says that $H_1(A_{n,k})=\Z_2$ when $n>2$. This agrees with the recently proved 
conjecture that 
the homology $H_i(\A_{n,k})$ is independent of $n$ and $k$ for $n>>i$ \cite{hw}.

To show that $L_{n,k}$ is contractible, we show that it is isomorphic to a subcomplex of $X_{n+3k}^\Theta$, the set of fixed points of $X_{n+3k}$ 
under the action of a certain finite subgroup $\Theta$ of $\Aut(F_{n+3k})$. 
Although $\A_{n,k}$ is 
not a subgroup of $\Aut(F_{n+k})$, it is a subgroup of the normalizer of $\Theta$ in $\Aut(F_{n+3k})$. Hence $\A_{n,k}$ acts on 
$X_{n+3k}^\Theta$, which is known to be contractible \cite{j1}. 
The space $L_{n+k}$ can be identified with the essential realization of an initial segment of 
the essential realization of $X_{n+3k}^\Theta$ under a certain norm, in the language of \cite{kv}.
Our  Theorem~\ref{thm:initseg}  says that the initial segments of the essential realization of 
$X_n^G$ are contractible for any finite group $G\leqslant \Aut(F_n)$.

\vs

Our motivation is the study of the stable automorphism group $\Aut_\infty$, the colimit of the inclusions $\Aut(F_n)\hookrightarrow \Aut(F_{n+1})$. 
The homology and homotopy type of the classifying space of the stable mapping class group $\Gamma_\infty$ has been determined recently \cite{mt,mw,t},  
making use of the mapping class groups $\Gamma_{g,n}$ of surfaces with any number of boundary components, and the fact that the homology 
groups $H_*(\Gamma_{g,n})$ are independent of $g$ and $n$ when $g>>*$ \cite{h}.  
Little is known about $B\!\Aut_\infty$. 
It was shown recently that the homology of $\A_{n,k}$ is also independent of the number of boundaries in a stable range \cite{hw}. 
This means in particular that we can work with $\A_{n,k}$ instead of $\Aut(F_n)$ if we are only interested in stable results. 
The boundaries of $\A_{n,k}$ allow us to define gluing operations analogous to the very useful gluing operations for the mapping class groups.
These  give a wealth of new operations on the automorphism of free groups and are used
in \cite{w} to show that the map 
on classifying spaces $B\Ga_\infty^+\to B\!\Aut_\infty^+$ is an infinite loop map after Quillen's plus-construction.

\vs

The paper is organized as follows. 
Section~\ref{Anksec} recalls the definition of the groups $\A_{n,k}$ and describes their relationship to mapping class groups of surfaces.
Section~\ref{nksec} defines the space $L_{n,k}$ and calculates a range for 
the virtual cohomological dimension of $\A_{n,k}$. 
In Section~\ref{fixedsection}, we prove the contractibility of initial segments in a general setting.
Section~\ref{contract} shows that $L_{n,k}$ is an initial segment of $X_{n+3k}^\Theta$.
In Section~\ref{extra}, we exhibit the sequence  $F_{n,k}\to L_{n,k}\to L_n^k$.  
Finally, Section~\ref{pressec} gives the presentations and calculates $H_1$.

\subsection*{Acknowledgment}

The authors would like to thank Benson Farb for helpful conversations 
and the referee for interesting suggestions. 
The work was conducted while the first author was partially supported by Louisiana Board of Regents Research Competitiveness 
Subprogram Contract LEQSF-RD-A-39 
and second was supported by a Marie Curie Fellowship of the European Community under contract number HPMF-CT-2002-01925.

\section{The groups {\rm A}$_{n,k}$}\label{Anksec}

\begin{figure}[ht]
\centering
\begin{picture}(0,0)%
\includegraphics{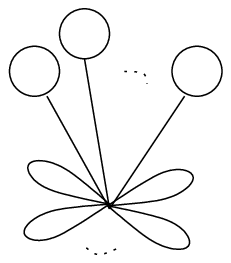}%
\end{picture}%
\setlength{\unitlength}{3158sp}%
\begingroup\makeatletter\ifx\SetFigFont\undefined%
\gdef\SetFigFont#1#2#3#4#5{%
  \reset@font\fontsize{#1}{#2pt}%
  \fontfamily{#3}\fontseries{#4}\fontshape{#5}%
  \selectfont}%
\fi\endgroup%
\begin{picture}(1425,1864)(676,-1319)
\put(826,-1261){\makebox(0,0)[rb]{\smash{\SetFigFont{10}{12.0}{\rmdefault}{\mddefault}{\updefault}{\color[rgb]{0,0,0}$x_2$}%
}}}
\put(1951,-811){\makebox(0,0)[lb]{\smash{\SetFigFont{10}{12.0}{\rmdefault}{\mddefault}{\updefault}{\color[rgb]{0,0,0}$x_{n}$}%
}}}
\put(2101, 14){\makebox(0,0)[lb]{\smash{\SetFigFont{10}{12.0}{\rmdefault}{\mddefault}{\updefault}{\color[rgb]{0,0,0}$y_k$}%
}}}
\put(676, 14){\makebox(0,0)[rb]{\smash{\SetFigFont{10}{12.0}{\rmdefault}{\mddefault}{\updefault}{\color[rgb]{0,0,0}$y_1$}%
}}}
\put(1201,389){\makebox(0,0)[rb]{\smash{\SetFigFont{10}{12.0}{\rmdefault}{\mddefault}{\updefault}{\color[rgb]{0,0,0}$y_2$}%
}}}
\put(751,-811){\makebox(0,0)[rb]{\smash{\SetFigFont{10}{12.0}{\rmdefault}{\mddefault}{\updefault}{\color[rgb]{0,0,0}$x_1$}%
}}}
\end{picture}
\caption{$G_{n,k}$}\label{Ggk}
\end{figure}

Let $G_{n,k}$ be the graph given in  Figure~\ref{Ggk}, obtained from a wedge of $n$ 
circles by attaching $k$ circles $C_1,\dots,C_k$ 
 on stems.  We think of the circles (together with the basepoint) as the boundary of the graph. 
The basepoint $*$ of $G_{n,k}$ is the central vertex.
Define 
 $$\A_{n,k}:=\pi_0\Htpy_*(G_{n,k};\del),$$
where $\Htpy_*(G_{n,k};\del)$ is the space of (basepointed) homotopy equivalences of
$G_{n,k}$ fixing the circles $C_j$ point-wise. By \cite[Prop. 0.19]{H},  homotopy equivalences fixing the circles are actually homotopy 
equivalences relative to the circles. This shows in particular that $A_{n,k}$ is a group. 

Labeling the first $n$ petals of $G_{n,k}$ with $x_1,\ldots, x_{n}$
and the $k$ circles with their stems $y_1, \ldots, y_k$ induces an 
identification of $\pi_1(G_{n+k})$ with  the free group $F_{n+k}:=$ $\langle x_1,\dots,x_n,y_1,\dots,y_k\rangle$. 
This leads to the following algebraic description of $\A_{n,k}$:
Define a map $\alpha:(F_{n+k})^{n+k} \to \operatorname{Hom}(F_{n+k}, F_{n+k})$ by 
$$\left\{\begin{tabular}{lcll}  $\alpha\langle\nu,w\rangle[x_i]=\nu_i$ & $1\le i\le n$\\
            $\alpha\langle\nu,w\rangle[y_j]=w_j^{-1}y_jw_j$ & $1\le j\le k$
\end{tabular}\right.$$
where 
$\langle \nu, w\rangle:= (\nu_1,\dots,\nu_{n},w_1,\dots,w_k)
\in (F_{n+k})^{n+k}$.

\begin{prop} There is a group isomorphism 
$$\A_{n,k}\cong\{\langle\nu,w\rangle
\in (F_{n+k})^{n+k}\ \big|\ 
\alpha\langle\nu,w\rangle\in \Aut(F_{n+k})\big\}$$
with the multiplication on the right hand side defined by  
$$\langle\nu,w\rangle . \langle\nu',w'\rangle=
   \langle\ \alpha\langle\nu',w'\rangle[\nu]\ ,\ w'.\,\alpha\langle\nu',w'\rangle[w]\ \rangle,$$
where we apply $\alpha\langle\nu',w'\rangle$ component-wise, and 
the unit is  $( x_1,\dots,x_n,1,\dots,1)$. 
\end{prop}

\begin{proof}
We use the labeling of $G_{n,k}$ by $x_i$'s and $y_j$'s as above (with a fixed orientation) to get an identification $\pi_1(G)\cong F_{n+k}$.
The ``lollies'' $y_j$ are made out of a stem, denoted $l_j$, from the basepoint to the circle, denoted $C_j$.

For any $f\in \Htpy_*(G;\del)$, we want to define an element 
$\langle\nu,w\rangle$ as above which depends only on the homotopy class 
of $f$. Set $\nu_i:=[f(x_i)]$ and  $w_j:=[l_{j}.f(\overline{l_j})]$ in $\pi_1(G)$, where $\overline{l_j}$ means the path $l_j$ 
taken in the reverse direction.  
Now $\al\langle\nu,w\rangle$ is the automorphism induced by $f$ on $\pi_1(G)$. 
Indeed, $[f(x_i)]=\nu_i$ and $[f(y_j)]=[f(l_j).f(C_j).f(\overline{l_j})]$ 
which is equal to $w_j^{-1}y_{j}w_j$ in $\pi_1(G)$. 
So $\langle\nu,w\rangle$ is an element of $\al^{-1}(\Aut(F_{n+k}))$ as required. 

Conversely, we want to show that any such $\langle\nu,w\rangle$ 
defines an element $[f]\in \pi_0\Htpy_*(G;\del)$. 
Define $f$ piecewise by $f(x_i)=\nu_i$, $f(l_j)=w_j^{-1}.l_{j}$ and 
$f(C_j)=C_{j}$. It is a continuous map, and a homotopy equivalence as   
it acts on $\pi_1(G)$ via $\al\langle\nu,w\rangle$. 

The two maps are homomorphisms  and inverse of one another.
Hence the two groups are isomorphic.
\end{proof}

The map $\alpha:\A_{n,k}\to \Aut(F_{n+k})$
is a homomorphism (where the multiplication of automorphisms $\phi.\psi$ is defined by 
first applying $\phi$ and 
then $\psi$.) Geometrically, it is the map that forgets that the circles were fixed.
It is not injective. There is in fact a short exact sequence
$$1\rar \Z^k \stackrel{i}{\rar} \A_{n,k} \stackrel{\alpha}{\rar} \A_n^k \to 1$$
where  $\A_n^k\cong\{\phi\in \Aut(F_{n+k})\ |\ \phi(y_j)\sim y_j \  \textrm{for}\ 1\le j\le k\} $ 
and where $i(z_1,\dots,z_k)=(x_1,\dots,x_n,y_1^{z_1},\dots,y_k^{z_k})$.  
Note that the group $\A_n^k$ can be described geometrically as $\pi_0\Htpy_*(G_{n,k};[\del])$, 
where $\Htpy_*(G_{n,k};[\del])$ is now the space of homotopy equivalences of $G_{n,k}$  fixing the circles up to a rotation.

The group $\Z^k$ maps to the center of $\A_{n,k}$.
These $k$ additional generators in $\A_{n,k}$ correspond to ``conjugating $y_i$ by itself''. 
Geometrically, these are Dehn twists along the $i$th boundary component in the corresponding mapping class 
group $\Ga_{g,k+1}$.

We will also consider the semidirect product  $A_{n,k}^\Si=\Si_k\ltimes \A_{n,k}$. 
Geometrically,  $\A_{n,k}^\Si:=\pi_0 \Htpy_*^\Si(G_{n,k};\del)$ with $\Htpy_*^\Si(G_{n,k};\del)$  the space of homotopy equivalences 
fixing the cycles $C_j$ up to a permutation. 
Algebraically, this means that the $y_i$'s can be permuted.

\vs

The groups $\A_{n,k}$ and $\A_n^k$ are closely related to mapping class groups when $n$ is even.
Indeed, note that the graph $G_{2g,k}$ is homotopy equivalent to a surface $S_{g,k+1}$ of 
genus $g$ with $k+1$ boundary circles and $\A_{n,k}$ can also be defined as the group of components of the homotopy equivalences of $S_{g,k+1}$ fixing 
the first $k$  boundary components point-wise and a point on the last one.
 
Let $\Ga_{g,k}:=\pi_0 \operatorname{Diff}^+(S_{g,k+1};\partial)$ denote the mapping class group of $S_{g,k+1}$  
with the boundaries point-wise fixed
and let $\Ga_g^{k+1}$ be the mapping class group of $S_{g,k+1}$ where the boundaries are fixed, but no longer point-wise. 
Fricke (1897) and Magnus (1934) for $g= 0,1$ and Nielsen (1927) for $g\ge 2$ proved that $\Gamma_g^{k+1}$ is the subgroup of 
$\operatorname{Out}(F_{2g+k})$ fixing the cyclic words defined by the boundaries \cite[p175]{m}.

Note that the inner automorphisms of $F_{n+k}$ form a normal subgroup of $\A_n^k$. Let 
$\operatorname{O}_n^k:=\A_n^k/\operatorname{Inn}(F_{n+k})$. So 
$\operatorname{O}_n^k=\{\phi\in \operatorname{Out}(F_{n+k})\ |\ \phi(y_j)\sim y_j \  \textrm{for}\ 1\le j\le k\}$. 
It follows that $\Ga_g^{k+1}$ is the subgroup of $\operatorname{O}_n^k$ fixing the  
conjugacy class of the last boundary component, that is the word $c:=[x_1,x_2]\dots [x_{2g-1},x_{2g}]\,y_1\dots y_k$. 
One actually has the following isomorphism of short exact sequences:
$$\xymatrix{ 1\ar[d] & 1\ar[d] \\
\Z^{k+1} \ar[d] \ar[r]^-{\cong} & \big\{(z,\phi)\in \Z^k\times \operatorname{Inn}(F_{2g+k})\ \big|\ \phi(c)=c\big\} \ar[d]\\
\Gamma_{g,k+1} \ar[r]^-{\cong} \ar[d] & \big\{\langle\nu,w\rangle\in \A_{2g,k}\ \big|\ \alpha\langle\nu,w\rangle[c]=c\big\} \ar[d] \\
\Gamma_g^{k+1} \ar[r]^-{\cong} \ar[d]  
 & \big\{\phi\in \operatorname{O}^k_{2g}\ \big|\ \phi(c)\sim c\big\} \ar[d]\\
1 & 1
}$$

\section{$(n,k)$-graphs}\label{nksec}

By a {\em graph}, we mean  a pointed, connected,
one-dimensional CW-complex.
Let $S^1$ be the circle with basepoint denoted by $\circ$.

\begin{defn} \label{defn:graph}
An {\em $(n,k)$-graph}  
$(\Gamma,\ul)$ is a graph $\Ga$ of genus $(n+k)$ with no separating edges, equipped with embeddings
 $l_j: S^1 \to \Gamma$ for $j=1,\ldots,k$ such that 
\begin{enumerate}
\item For every $j=1,\ldots,k$, the point $\circ_j:=l_j(\circ)$ is a vertex of $\Gamma$.
\item The vertices of $\Gamma$ have valence at least $2$
and if a vertex $v$ has valence $2$, then  $v=\circ_j$ for some $j$
or  $v$ is the basepoint $*$.
\item Any two $\SP_j:=l_j(S^1)$ meet in at
most a point and the dual graph of the $\SP_j$'s is a forest. 
\end{enumerate}
\end{defn}

\begin{figure}[ht]
\centering
\begin{picture}(0,0)%
\includegraphics{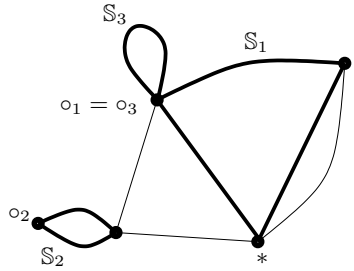}%
\end{picture}%
\setlength{\unitlength}{3108sp}%
\begingroup\makeatletter\ifx\SetFigFont\undefined%
\gdef\SetFigFont#1#2#3#4#5{%
  \reset@font\fontsize{#1}{#2pt}%
  \fontfamily{#3}\fontseries{#4}\fontshape{#5}%
  \selectfont}%
\fi\endgroup%
\begin{picture}(2566,2149)(923,-1784)
\put(946,-1366){\makebox(0,0)[rb]{\smash{\SetFigFont{9}{10.8}{\rmdefault}{\mddefault}{\updefault}{\color[rgb]{0,0,0}$\circ_2$}%
}}}
\put(1801,-511){\makebox(0,0)[rb]{\smash{\SetFigFont{9}{10.8}{\rmdefault}{\mddefault}{\updefault}{\color[rgb]{0,0,0}$\circ_1=\circ_3$}%
}}}
\put(2836,-1726){\makebox(0,0)[rb]{\smash{\SetFigFont{9}{10.8}{\rmdefault}{\mddefault}{\updefault}{\color[rgb]{0,0,0}$*$}%
}}}
\put(2611,-16){\makebox(0,0)[lb]{\smash{\SetFigFont{9}{10.8}{\rmdefault}{\mddefault}{\updefault}{\color[rgb]{0,0,0}$\SP_1$}%
}}}
\put(1486,209){\makebox(0,0)[lb]{\smash{\SetFigFont{9}{10.8}{\rmdefault}{\mddefault}{\updefault}{\color[rgb]{0,0,0}$\SP_3$}%
}}}
\put(991,-1726){\makebox(0,0)[lb]{\smash{\SetFigFont{9}{10.8}{\rmdefault}{\mddefault}{\updefault}{\color[rgb]{0,0,0}$\SP_2$}%
}}}
\end{picture}
\caption{Example of a $(2,3)$-graph}\label{nkgraph}
\end{figure}

We call the embedded circle $(\SP_j,\circ_j)$ a {\em cycle}. 
In this way, an $(n,k)$-graph is a graph of genus 
$n+k$ with $k$ cycles. (See Figure~\ref{nkgraph} for an example.)

Let $(G_{n,k},\underline{\kappa})$ be the generalized 
$(n,k)$-graph, where
$G_{n,k}$ was defined in Section~\ref{Anksec} (see Figure~\ref{Ggk}) and 
$\kappa_j:S^1\to G_{n,k}$ identifies $S^1$ with the $j$th circle $C_j$. 
Note that $G_{n,k}$
has $k$ separating edges, and hence this is not an actual 
$(n,k)$-graph.

\begin{defn}
A {\em marked $(n,k)$-graph} $(\Ga,\phi)$ is a graph $\Ga$ equipped with a pointed 
homotopy equivalence $\phi:G_{n,k} \to \Ga$ 
such that 
$(\Ga,\phi\circ \underline{\kappa})$ is a $(n,k)$-graph.

Two marked $(n,k)$-graphs $(\Gamma^1,\phi^1)$ and $(\Gamma^2,\phi^2)$ are
{\em equivalent} if there is a graph isomorphism 
$\psi: \Gamma^1 \to \Gamma^2$ 
such that $\psi\circ \phi^1(C_j)=\phi^2(C_j)$ 
and $\psi\circ \phi^1 \simeq \phi^2$ (rel. to $C_j$ for all $j$.) 
\end{defn}

The intuition (which is made precise in the proof of Theorem \ref{thm:posetequivalence}) 
is that 
condition (3) in Definition~\ref{defn:graph}  prevents the existence of
Whitehead moves (blowing an edge and then collapsing it) which would take a 
`good marking' to a `bad marking'.

For an $(n,k)$-graph $(\Gamma,\ul)$, let $\hat \Gamma$ be 
the graph of genus $n$ obtained from $\Ga$ by collapsing
all the edges in the cycles $\SP_1,\dots,\SP_k$.
Let $e$ be an edge
of $\Gamma$ 
which does not define a loop in $\Ga$ or in $\hat\Gamma$. 
Then composition with the edge collapse  $\Ga \to \Gamma/e$ induces embeddings
$l_j/e:S^1\to\Ga/e$ such that  $(\Ga/e,\ul/e)$ is again a $(n,k)$-graph. 
By an edge collapse in an $(n,k)$-graph, we will always mean the
collapse of an edge satisfying the above hypothesis.  
For a marked $(n,k)$-graph, the marking of the collapsed graph is also obtained by composition with the collapse.

A  forest collapse in a $(n,k)$-graph $(\Ga,\ul)$ 
is a sequence of edge collapses. As the order of the collapses does not matter, it is defined by the union of the edges to 
be collapsed, which is a forest $F$ in $\Ga$.  
We denote the collapsed graph by $(\Gamma/F,\ul/F)$.

\vs

Let $L_{n,k}$ denote the set of equivalence classes of 
marked $(n,k)$-graphs. 
Define a poset structure on $L_{n,k}$ by saying that $(\Ga^1,\phi^1)\le (\Ga^2,\phi^2)$ if there is 
a forest  $F$ in $\Ga^2$ such that $(\Ga^2/F,\phi^2/F)$ is equivalent to $(\Ga^1,\phi^1)$. 
By abuse of notation, we will also denote by $L_{n,k}$ 
 the simplicial complex which is the geometric realization of this poset.

Using the geometric description of $A^\Si_{n,k}$, we obtain a
natural action of $A^\Si_{n,k}$ (and thus of $\A_{n,k}$) on $L_{n,k}$: 
For $\lambda\in A^\Si_{n,k}$
and $(\Ga,\phi)\in L_{n,k}$, define 
$$\lambda\cdot(\Ga,\phi):=(\Ga,\phi\circ\lambda^{-1}).$$
We consider two equivalence relations on $(n,k)$-graphs:
Two $(n,k)$-graphs $(\Gamma^1,\ul^1)$ and $(\Gamma^2,\ul^2)$ 
are {\em $\Si$-equivalent} if there is a graph isomorphism
$\phi: \Gamma^1 \to \Gamma^2$ and a permutation
$\sigma\in\Sigma_k$
such that  $\phi\circ l_j^1(\circ)= l_{\sigma(j)}^2(\circ)$
and $\phi\circ l_j^1\simeq l_{\sigma(j)}^2$ (rel. to $\circ$) for each
$j=1,\ldots,k$. 
The $(n,k)$-graphs are {\em equivalent} if the permutation $\sigma$ is trivial.

Let $QL_{n,k}^\Si$ and $QL_{n,k}$ denote the set of all equivalence
classes of $(n,k)$-graphs under each equivalence relation. 
Define spaces $QL_{n,k}^\Si$ and $QL_{n,k}$ by taking the elements of the sets as $0$-simplices and 
attaching a $p$-simplex for each equivalence class of chain of forests 
$\emptyset \not = F_1 \subsetneq F_2 \subsetneq \cdots \subsetneq F_p$,
where two chains $F^1_i$ and $F^2_i$ are equivalent  if there is a ($\Si$-)equivalence
$\phi: (\Gamma,\ul) \to (\Gamma,\ul)$ such that $\phi(F^1_i)=F^2_i$
for all $i=1,\ldots, p$. 
Note that these spaces are not simplicial complexes in general.
 The space $QL_{n,k}$ is the 
moduli space of labeled $(n,k)$-graphs, that is $(n,k)$-graphs with cycles labeled $1,\dots,k$, whereas in $QL^\Si_{n,k}$ the 
cycles are unlabeled.

\begin{prop}
The action of $A^\Si_{n,k}$ (resp.\ $\A_{n,k}$) on $L_{n,k}$ is transitive on the sets of marked $(n,k)$-graphs
with the same underlying unlabeled 
$(n,k)$-graph (resp.~same underlying labeled $(n,k)$-graph), 
with finite simplex stabilizers, and  the quotient of the action  is $QL^\Si_{n,k}$ (resp.\ $QL_{n,k}$).
\end{prop}

We will prove the following theorem in the next
two sections.

\begin{thm} \label{thm:main}
The space $L_{n,k}$ is contractible.
\end{thm}
\eject

We use this result to give a range for the virtual cohomological dimension of $\A_{n,k}$. 
We recall that it is already known that 
the vcd of $\Aut(F_n)=\A_{n,0}$ is $2n-2$ \cite{cv} and that the
vcd of $\A_0^k$ is $k-1$ \cite{c}. 
Since $A_{n,k}$ is a subgroup of $\Aut(F_{n+3k})$ (see 
Proposition \ref{prop:normalizerinclusion})
and the latter group is
virtually torsion free and of finite vcd, the former group is as well
\cite[Example~1 on page 229]{b}.

\begin{thm}\label{vcd}
\begin{enumerate}
\item $2n+2k-1\le vcd(\A_{n,k})\le 2n+3k-2\ \ $ if $k \geq 1$.
\item $vcd(\A_{n,1})=2n+1$.
\item $vcd(\A_{0,k})=2k-1$.
\end{enumerate}
\end{thm}

\begin{proof}[Proof of Theorem \ref{vcd}] 
For (1), we first note that there is a copy of $\Z^{2n+2k-1}$ in $\A_{n,k}$:
take $n$ generators sending
$x_i$ to $x_i y_k$ and fixing all other generators, then $n$ generators 
sending  $x_i$ to $y_k x_i$, then $k$ conjugating $y_j$ by $y_k$ for $1\le j\le k$
and finally $k-1$ conjugating $y_j$ by $y_j$ for $j \not = k$. 

For an upper bound, we calculate the dimension of a maximal simplex in 
$QL_{n,k}$.  Suppose the cycles lie in the maximally blown up
graph in $c$ connected components and that the graph has $v$ vertices. 
The $k+1$ vertices $*, \circ_1, \ldots, \circ_k$ all have valence 2.  Since
two cycles must meet in a valence 4 vertex and the dual graph of the cycles in the 
graph is a forest, there are $k-c$ vertices of valence 4 in the graph.  The remaining
vertices are valence 3.  Hence there are $e=3(v-2k+c-1)/2+(k+1)+2(k-c)=(3v-c-1)/2$
edges.  Since $v-e=1-(n+k)$, this yields $v=2n+2k+c-1$ vertices.  We can collapse
these $v$ vertices down to 1 by first collapsing all except one edge in each cycle and
then collapsing some remaining maximal tree.  This yields a simplex of dimension
$2n+2k+c-2$.  Since $c \leq k$, the dimension of a maximal simplex is less
than or equal to $2n+3k-2$. The result follows from Theorem \ref{thm:main}.

For (2), let $k=1$ in (1). 

Finally for (3) if $n=0$ then $c=1$ and so $v=2n+2k+c-1=2k$.  Hence a
maximal simplex has dimension $2k-1$. 
\end{proof}

\begin{figure}[ht]
\centerline{\begin{picture}(0,0)%
\includegraphics{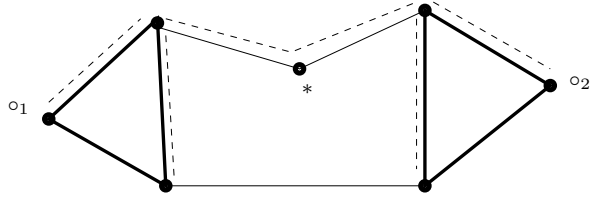}%
\end{picture}%
\setlength{\unitlength}{2763sp}%
\begingroup\makeatletter\ifx\SetFigFont\undefined%
\gdef\SetFigFont#1#2#3#4#5{%
  \reset@font\fontsize{#1}{#2pt}%
  \fontfamily{#3}\fontseries{#4}\fontshape{#5}%
  \selectfont}%
\fi\endgroup%
\begin{picture}(4800,1729)(451,-1928)
\put(5251,-961){\makebox(0,0)[lb]{\smash{\SetFigFont{8}{9.6}{\rmdefault}{\mddefault}{\updefault}{\color[rgb]{0,0,0}$\circ_2$}%
}}}
\put(451,-1186){\makebox(0,0)[rb]{\smash{\SetFigFont{8}{9.6}{\rmdefault}{\mddefault}{\updefault}{\color[rgb]{0,0,0}$\circ_1$}%
}}}
\put(2851,-1036){\makebox(0,0)[lb]{\smash{\SetFigFont{8}{9.6}{\rmdefault}{\mddefault}{\updefault}{\color[rgb]{0,0,0}$*$}%
}}}
\end{picture}}
\caption{\label{lollyfig2}
A maximal simplex in $QL_{1,2}$ with a 6 edge forest}
\end{figure}

 In Figure \ref{lollyfig2},
we show an example of a maximal simplex of $QL_{1,2}$ given by a graph with 2 connected cycle components 
and a forest consisting of $6>5$ edges. 

\section{Initial segments of fixed point subcomplexes}\label{fixedsection}

Recall from \cite{hv} 
that the spine of Auter space $X_n$ is a poset 
of pointed marked graphs of genus $n$, where the marking is given by a map from the rose
$R_n:=\vee_n S^1$. $X_n$ is contractible and admits an action of $\Aut(F_n)$. 
Let $G$ be a finite subgroup of $\Aut(F_n)$ and consider the fixed point 
subcomplex $X_n^G$. A point in 
$X_n^G$ is a marked graph with a $G$-action.
$X_n^G$ is contractible \cite[Thm.~1.1]{j1}.

Edge collapses induce a poset structure on $X^G_n$ with minimal elements 
the {\em reduced} marked
graphs, i.e. those with no $G$-invariant
subforest.  An edge $e$ of a marked $G$-graph $\Gamma$ is
{\em essential} if there exists a maximal $G$-invariant forest which does not 
contain $e$. 
The graph $\Gamma$ is {\em essential}
if all of its edges are essential.  Let $E_n^G$ 
be the full subcomplex
of $X_n^G$ spanned by essential marked graphs.  
Collapsing inessential edges induces a poset map $f:X_n^G\to E_n^G$ such that 
$f(x)\leq x$ for all $x$ in $X_n^G$, and hence a
deformation retraction on the corresponding spaces by the Poset Lemma 
(see \cite{q} or \cite[Lem.~4.1]{c} for a direct proof of this special case).

Let $\mathcal{W}$ be the set of conjugacy classes of elements of $F_n$. 
Recall from \cite[Section~2]{j1} 
that, given any well ordering of $\mathcal{W}$ and of $F_n$, 
we can define a norm
$$\|\cdot\|_{tot}=\|\cdot\|_{out} \times \|\cdot\|_{aut}
\in \Z^{\mathcal{W}} \times \Z^{F_n}$$
on $E_n^G$ 
which well orders the reduced marked graphs of $E_n^G$ (using the corresponding 
lexicographic ordering). 
For $w \in \mathcal{W}$, the $w$th
component of the norm is the sum of the cyclic 
lengths of the paths $\phi(xw)$
for all $x \in G$, where the cyclic length is the number of edges in the
path reduced cyclically, i.e. the length of the shortest path in its unpointed 
homotopy class.  
And for $z \in F_n$, the $z$th
component of the norm is the sum of the lengths
of the paths $\phi(xz)$ for all $x \in G$, where we now consider the
shortest paths in the pointed homotopy classes.

\vs

A nonempty set $\Lambda$ of reduced marked
graphs of $E_n^G$ is called an {\em initial segment} if
whenever two reduced marked graphs $r_1, r_2$ satisfy
$\|r_1\|_{tot} \leq \|r_2\|_{tot}$ and $r_2 \in \Lambda$,
then $r_1 \in \Lambda$.  The {\em realization} of an initial
segment $\Lambda$ is the space 
$$\|\Lambda\|=\bigcup_{r\in\Lambda}st(r),$$ 
the union of the ascending stars (in $E_n^G$) of the reduced
marked graphs in $\Lambda$.

Let $\Lambda$ be an initial segment and let $\Ga$ be a marked graph
in $\|\Lambda\|$. An edge $e$ of $\Gamma$
 is {\em essential in $\Lambda$} if
there is a maximal $G$-invariant subforest $\mathcal{F}$
of $\Gamma$ such that $e \not \in \mathcal{F}$ and
$\Gamma/\mathcal{F} \in \Lambda$.
The graph $\Gamma$ is {\em essential in $\Lambda$}
if all of its edges are essential in $\Lambda$.
The {\em essential realization} of an initial
segment $\Lambda$ is the full subcomplex $|\Lambda|$
of $\|\Lambda\|$ spanned by the essential marked graphs 
of $\Lambda$.

\vs

Let $\Gamma$ be a 
marked graph in $E_n^G$
and let $v$ be a vertex of $\Gamma$.
An {\em ideal edge} $\ga$  at the vertex $v$ is 
a set of half edges of $\Ga$ terminating at $v$ such that 
the {\em blow-up} $\Ga^\ga$ is again in $E_n^G$, 
where  $\Gamma^{\gamma}$  is the
graph obtained by $G$-equivariantly pulling the
half edges in $\gamma$ away from $v$
(see \cite[Section~5.A.]{kv}  and \cite[Section~2]{j1} for a characterization of 
ideal edges when $\Ga$ is reduced.) 
More precisely, one creates 
a new edge orbit $G\gamma$ and a new vertex orbit
$Gv(\gamma)$.  In $\Gamma^{\gamma}$,
$\gamma$ is an actual edge that goes from $v(\gamma)$
to $v$.  In addition, each half edge $e \in \gamma$ 
is now attached to
$v(\gamma)$ instead of $v$, and this process is done equivariantly on
$G\gamma$.   
Note that the original graph $\Ga$ can be recovered
from its blow up $\Ga^\gamma$ by collapsing $G\gamma$.

An {\em ideal forest} in a reduced 
marked graph is a sequence of ideal edges satisfying 
certain technical conditions, which make sure that by blowing up these ideal
edges, one can construct a blow-up $\Ga^\mathcal{F}$ in the star of $\Ga$ 
with an actual $G$-equivariant forest, whose
partial collapses contain all the blow-ups of the ideal edges in the ideal forest.
One can define a poset structure on the set of ideal forests of a reduced
marked graph such that blowing up induces an isomorphism between this poset and 
the star of $\Gamma$ in $E_n^G$ \cite[Prop.~5.9]{kv}.

Let $\ga$ be an ideal edge of a reduced 
$G$-graph $\Ga$. Define $D(\ga)$ to be the set of edges $c$
of $\ga$ such that $Gc$ is a forest in $\Ga^\ga$.
To a pair  $(\gamma,c)$ with $c\in D(\ga)$ corresponds 
a {\em Whitehead move} that consists of first blowing up
$\gamma$ and then collapsing $Gc$. 
An ideal edge $\gamma$ is
{\em reductive} if there is a $c\in D(\ga)$ such that 
performing the Whitehead move $(\gamma,c)$
yields a graph of smaller norm.

Let $Y$ be any of $E_n^G$, $\|\Lambda\|$ or $|\Lambda|$. 
The {\em reductive link} of a reduced marked graph $\Gamma$ in $Y$ 
is the part of the star of $\Gamma$ in $Y$ spanned by
ideal forests with at least one reductive ideal edge.
The {\em purely reductive link} of $\Gamma$ in $Y$ is the
part of the star of $\Gamma$ in $Y$ spanned by nontrivial
ideal forests with all edges reductive. By the Poset Lemma,
there is a deformation
retraction from the reductive link of $\Gamma$ to the purely
reductive link of $\Gamma$ given by collapsing non-reductive
ideal edges.

\begin{thm} \label{thm:initseg} 
If $\Lambda$ is an initial segment of $E_n^G$, then
both its realization $\|\Lambda\|$ and its essential realization 
$|\Lambda|$ are contractible.
\end{thm}

\begin{proof}
By the Poset Lemma, collapsing inessential edges induces a 
deformation retraction from $\|\Lambda\|$ to $|\Lambda|$.   
Thus it is enough to show that $\|\Lambda\|$ is 
contractible.

If $r$ is a reduced marked graph in $\|\Lambda\|$, we have already
observed that the reductive link of $r$ in $\|\Lambda\|$
deformation retracts onto the purely reductive link of $r$
in $\|\Lambda\|$.  Since $\Lambda$ is an initial segment, the purely
reductive link of $r$ in $\|\Lambda\|$ is the purely reductive
link of $r$ in $E_n^G$.  From \cite[Lem.~3.1 - 3.7]{j1},
 this reductive link (referred to as $S(\mathcal{R})$ in \cite{j1}) is contractible.
As $\Lambda$ is an initial segment, all $\Ga'$ such that $\|r'\|<\|r\|$ are
also in $\Lambda$.
The reductive link of a reduced marked graph $r$ is the
intersection
$$red(r)=st_{\|\Lambda\|}(r)\cap\bigcap_{\|r'\|<\|r\|}st_{\|\Lambda\|}(r').$$
Now the norm well orders the reduced marked graphs in
$\|\Lambda\|$. We can build $\|\Lambda\|$ as the union of the stars of the
elements of $\Lambda$, starting with the star of the smallest element and
adding the stars successively. 
As the star of a graph is contractible and all the successive intersections are
contractible, 
transfinite induction implies that every connected
component of $\|\Lambda\|$ is contractible.

It remains to show that $\|\Lambda\|$ is connected. 
From \cite[Thm.~2]{k}, we know that any two reduced marked graphs
in $E_n^G$ are connected by a sequence of Whitehead moves.
The reduced marked graphs of $E_n^G$ are ordered by the norm. 
Choose $I \subseteq \Z^{\mathcal{W}} \times \Z^{F_n}$ such that  
the ordered set of reduced marked graphs of $E_n^G$ is $\{r_i : i \in I\}$ 
and denote by $r_0$ the least element. 
Note that $r_0 \in \Lambda$ because $\Lambda$ is nonempty.  
Note also that the stars $st_{\|\Lambda\|}(r_i)$ and
 $st_{E_n^G}(r_i)$ of $r_i$ in $\|\Lambda\|$ and $E_n^G$ are equal. 
We will denote these stars by $st(r_i)$.

Suppose that $\|\Lambda\|$ is not connected.
The star $st(r_0)$ is contractible and
thus connected. Take the least index $j$ such that
$\cup_{i \leq j} st(r_i)$ is not connected.
Then we can get from $r_0$ to any
$r_i$ with $i<j$ in $\|\Lambda\|$, but not from $r_0$
to $r_j$.
Since $E_n^G$ is connected and $I$ is a well ordering, we can choose a
least $k$ such that there is a path from $r_0$ to $r_j$ in
$\cup_{i \leq k} st(r_i)$.
This implies that the reductive link of $r_k$ is not
connected, however, since it has at least two components: one coming from
$r_0$ and the other coming from $r_j$.  This contradicts the
fact that the reductive link is contractible.
So $\|\Lambda\|$ is connected.
\end{proof}

\section{Contractibility of $L_{n,k}$}\label{contract}

In this section, we prove Theorem \ref{thm:main} by identifying
$L_{n,k}$ with the essential realization of an initial
segment and applying Theorem \ref{thm:initseg}.

\vs

\begin{figure}[ht]
\centering
\begin{picture}(0,0)%
\includegraphics{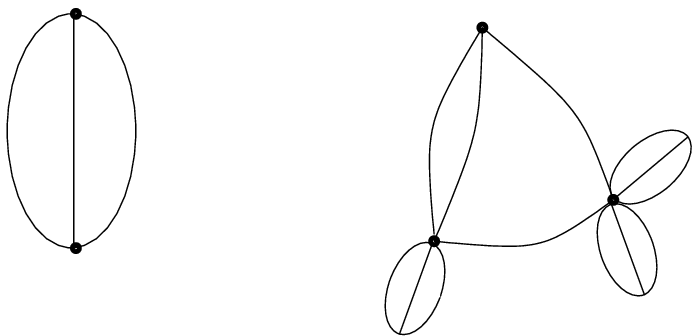}%
\end{picture}%
\setlength{\unitlength}{2901sp}%
\begingroup\makeatletter\ifx\SetFigFont\undefined%
\gdef\SetFigFont#1#2#3#4#5{%
  \reset@font\fontsize{#1}{#2pt}%
  \fontfamily{#3}\fontseries{#4}\fontshape{#5}%
  \selectfont}%
\fi\endgroup%
\begin{picture}(4523,2396)(451,-1761)
\put(991,-466){\makebox(0,0)[lb]{\smash{\SetFigFont{9}{10.8}{\rmdefault}{\mddefault}{\updefault}{\color[rgb]{0,0,0}$e^j_2$}%
}}}
\put(451,-241){\makebox(0,0)[rb]{\smash{\SetFigFont{9}{10.8}{\rmdefault}{\mddefault}{\updefault}{\color[rgb]{0,0,0}$e^j_1$}%
}}}
\put(1396,-241){\makebox(0,0)[lb]{\smash{\SetFigFont{9}{10.8}{\rmdefault}{\mddefault}{\updefault}{\color[rgb]{0,0,0}$e^j_3$}%
}}}
\put(4366,-826){\makebox(0,0)[rb]{\smash{\SetFigFont{9}{10.8}{\rmdefault}{\mddefault}{\updefault}{\color[rgb]{0,0,0}$\circ_1=\circ_3$}%
}}}
\put(3241,-1051){\makebox(0,0)[rb]{\smash{\SetFigFont{9}{10.8}{\rmdefault}{\mddefault}{\updefault}{\color[rgb]{0,0,0}$\circ_2$}%
}}}
\put(946,-1366){\makebox(0,0)[lb]{\smash{\SetFigFont{9}{10.8}{\rmdefault}{\mddefault}{\updefault}{\color[rgb]{0,0,0}$\circ_j$}%
}}}
\put(901,479){\makebox(0,0)[lb]{\smash{\SetFigFont{9}{10.8}{\rmdefault}{\mddefault}{\updefault}{\color[rgb]{0,0,0}$v_j$}%
}}}
\end{picture}
\caption{The graph $\theta^j$ and a graph of $Fix$}\label{Fixfig}
\end{figure}

Let $m=n+3k$.  Let
$\{x_1,\dots,x_{n},y_1,\dots,y_k,u_1,\dots,u_k,v_1,\dots,v_k\}$ denote the
generators of $F_m$. 
For each $j=1,\dots,k$, let $\Theta_j\cong \Z/3$ be the subgroup of $\Aut(F_m)$ 
generated by $\theta_j$, where  
$\theta_j(u_j)=v_j$, $\theta_j(v_j)=v_j^{-1} u_j^{-1}$, and 
all other generators are fixed. 
This group can be realized as the group of rotations of
the three edges of a theta-graph (see Figure~\ref{Fixfig}). 
Let
$$\Theta = \Theta_1 \times \Theta_2 \times \cdots \times \Theta_k.$$
be the direct product of these $k$ groups.
The normalizer $N_{\Aut(F_m)}(\Theta)$
acts on the fixed point subcomplex $X_m^{\Theta}$ and thus on $E_m^{\Theta}$. Both
spaces are contractible \cite{j1} (see also Section \ref{fixedsection}).  

Note that a similar fixed point subcomplex was used in \cite{j2} to study the holomorph of free groups $F_n \rtimes \Aut(F_n)$.
This last group is the subgroup of $\A_n^1$ of automorphisms where $y_1$ does not occur in the images of the $x_i$'s. 
One could work with circles and a $\Z/2$-action, instead of $\Theta$-graphs and a $\Z/3$-action. 
This would however bring a few extra technical 
problems (having to introduce extra vertices to avoid acting on edges 
by inversions) while producing the same complex. (See also Section~\ref{extra}.)

Although $\A_{n,k}$ is not a subgroup of $\Aut(F_{n+k})$, it is a subgroup of 
$\Aut(F_{n+3k})$:

\begin{prop}\label{prop:normalizerinclusion}
There is an inclusion of groups 
$\beta:A^\Si_{n,k}\hookrightarrow N_{\Aut(F_{n+3k})}(\Theta)\leqslant \Aut(F_{n+3k})$.
\end{prop}

\begin{figure}[ht]
\centering
\begin{picture}(0,0)%
\includegraphics{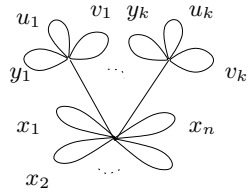}%
\end{picture}%
\setlength{\unitlength}{2960sp}%
\begingroup\makeatletter\ifx\SetFigFont\undefined%
\gdef\SetFigFont#1#2#3#4#5{%
  \reset@font\fontsize{#1}{#2pt}%
  \fontfamily{#3}\fontseries{#4}\fontshape{#5}%
  \selectfont}%
\fi\endgroup%
\begin{picture}(1675,1639)(576,-1319)
\put(826,-1261){\makebox(0,0)[rb]{\smash{\SetFigFont{9}{10.8}{\rmdefault}{\mddefault}{\updefault}{\color[rgb]{0,0,0}$x_2$}%
}}}
\put(1951,-811){\makebox(0,0)[lb]{\smash{\SetFigFont{9}{10.8}{\rmdefault}{\mddefault}{\updefault}{\color[rgb]{0,0,0}$x_{n}$}%
}}}
\put(751,-811){\makebox(0,0)[rb]{\smash{\SetFigFont{9}{10.8}{\rmdefault}{\mddefault}{\updefault}{\color[rgb]{0,0,0}$x_1$}%
}}}
\put(1951,164){\makebox(0,0)[lb]{\smash{\SetFigFont{9}{10.8}{\rmdefault}{\mddefault}{\updefault}{\color[rgb]{0,0,0}$u_k$}%
}}}
\put(2251,-361){\makebox(0,0)[lb]{\smash{\SetFigFont{9}{10.8}{\rmdefault}{\mddefault}{\updefault}{\color[rgb]{0,0,0}$v_k$}%
}}}
\put(751, 89){\makebox(0,0)[rb]{\smash{\SetFigFont{9}{10.8}{\rmdefault}{\mddefault}{\updefault}{\color[rgb]{0,0,0}$u_1$}%
}}}
\put(676,-361){\makebox(0,0)[rb]{\smash{\SetFigFont{9}{10.8}{\rmdefault}{\mddefault}{\updefault}{\color[rgb]{0,0,0}$y_1$}%
}}}
\put(1351,164){\makebox(0,0)[rb]{\smash{\SetFigFont{9}{10.8}{\rmdefault}{\mddefault}{\updefault}{\color[rgb]{0,0,0}$v_1$}%
}}}
\put(1426,164){\makebox(0,0)[lb]{\smash{\SetFigFont{9}{10.8}{\rmdefault}{\mddefault}{\updefault}{\color[rgb]{0,0,0}$y_k$}%
}}}
\end{picture}
\caption{Graph giving the inclusion $\A_{n,k}\hookrightarrow \Aut(F_{n+3k})$}\label{G3}
\end{figure}

\begin{proof}
The map $\beta$ is induced geometrically
by attaching two additional circles at the end of each of the $k$ stems (see Figure~\ref{G3}). 
Algebraically, this map can be described as follows:
for $\langle \s,\nu_i,w_j\rangle\in A^\Si_{n,k}=\Si_k\rtimes A_{n,k}$, 
$$\phi=i\langle \s,\nu_i,w_j\rangle\ \ \textrm{maps}\ \  \left\{\begin{tabular}{lcl}
$x_i$ & $\rar$ & $\nu_i$ \\
$y_j$ & $\rar$ & $w_j^{-1}\ y_{\sigma(j)}\ w_j$ \\
$u_j$ & $\rar$ & $w_j^{-1}\ u_{\sigma(j)}\ w_j$ \\
$v_j$ & $\rar$ & $w_j^{-1}\ v_{\sigma(j)}\ w_j$
\end{tabular}\right.$$
One then checks easily that $\phi=\beta\langle \s,\nu_i,w_j\rangle$
satisfies 
 $\phi\ \theta_l\ \phi^{-1} = \theta_{\sigma(l)}$
for each $l = 1, \ldots, k$. 
\end{proof}

For each $j=1,\ldots,k$, let $\theta^j$ denote a graph with
two vertices $\circ_j$ and $v_j$, and three edges
$e^j_1, e^j_2, e^j_3$ where each $e^j_i$ goes from $\circ_j$ to
$v_j$ (See Figure~\ref{Fixfig}).
 Consider the subcomplex $Fix$ of $X_m$ spanned by marked graphs
$\phi: R_n \to \tilde \Gamma$ which satisfy the following conditions:
\begin{itemize}
\item The graph $\tilde \Gamma$ is obtained from a graph
$\Gamma$ of genus $n+k$ 
by attaching the graphs $\theta^1, \ldots, \theta^k$ to $\Gamma$, 
identifying the vertex $\circ_j$ of $\theta^j$ with a vertex of $\Ga$ and
remembering these ---possibly multiple--- labels.
\item For each $j=1,\ldots,k$, there exists a path $w_j$ in
$\tilde \Gamma$ from $\circ_j$ to $*$ such that
$$\phi(u_j) = \bar w_j \star e^j_1 \star \bar e^j_2 \star  w_j$$
and
$$\phi(v_j) = \bar w_j \star e^j_2 \star \bar e^j_3 \star  w_j$$
where $\bar p$ denote the path $p$ with the reverse orientation and $\star$
denotes the concatenation of paths. 
\end{itemize}

\begin{prop} $E^\Theta_m \cong Fix \cong X^\Theta_m.$
\end{prop}

\begin{proof}
Each marked graph in $Fix$ is fixed by the
action of $\Theta$ as $\theta_j$ cyclically
rotates the edges of the attached $\theta^j$. So $Fix\subset X^\Theta_m$.

Let $R^\Theta_{n,k}$ be the graph obtained by wedging the $k$ graphs 
$\theta^j$ to the rose $R_{n+k}$. 
Recall that a graph is reduced if it has no invariant subforest. Note that all the
reduced graphs in $Fix$ have the same underlying graph $R^\Theta_{n,k}$.

Krsti\'c showed that any two reduced $\Theta$-graphs in
$X_m^{\Theta}$ are connected by a sequence of
Whitehead moves $(\{e,f\}, f)$ where $e,f$ have the
same terminal vertex and $Stab_{\Theta}(e) \subseteq
Stab_{\Theta}(f)$  (and $e\notin \Theta f$)
\cite[Prop.~4']{k}.  Pictorially
the move $(\{e,f\}, f)$ corresponds
to pulling each edge in $\Theta e$ along the corresponding
edge in $\Theta f$.  Because of the condition on stabilizers,
we can only perform Whitehead moves on $R_{n,k}^\Theta$
of the form $(\{e,f\}, f)$ or $(\{\bar e^j_i,f\}, f)$
where $e,f$ are edges of $R_{n+k}$ and $\bar e^j_i$ is in
$\theta^j$. The first type of moves affect only the first $n+k$ petals,
whereas the second type takes the subgraph $\theta^j$ around the edge $f$.
The reduced graphs in $Fix$ are closed under these moves.

Hence the reduced graphs in $X_m^{\Theta}$ are exactly the
reduced graphs in $Fix$.  Now the ascending star of a reduced
graph in $Fix$ is the same as the ascending star of the
graph in $X_m^{\Theta}$ as equivariant blow-ups preserve the $\theta$-graphs 
(see \cite[Claim 6.6]{j0} for a similar computation). 
So $Fix=X_m^\Theta$. 
Note finally that all graphs in $Fix$ are essential 
as no edge
in a $\theta^j$ can be collapsed $\Theta$-equivariantly and the other edges,
being non-separating by hypothesis, are thus also essential. 
So $Fix=X^\Theta_m=E^\Theta_m$.
\end{proof}

We now choose a well ordering of $\Z^{\mathcal{W}} \times \Z^{F_n}$ such that
the the space $L_{n,k}$ of marked $(n,k)$-graphs 
is the realization of an
initial segment of the induced norm $\|\cdot\|_{tot}$.
Well order the set of conjugacy classes $\mathcal{W}$ so that its first $k$ elements are
$$u_1 y_1 u_1, u_2 y_2 u_2, \ldots, u_k y_k u_k$$
and well order $F_n$ in any manner.
Now define $\Lambda$ to be the set of 
reduced $\Theta$-graphs in $E^\Theta_m$ whose norm
is minimal in the first $k$ coordinates. This is an initial segment in the
lexicographic ordering induced by the norm.
Recall that all reduced marked graphs in
$E^\Theta_m$ have underlying graph $R^\Theta_{n,k}$. 
The minimal norm on the component $u_j y_j u_j$ is thus 
$(2+1+2) \cdot | \Theta | = 5 \cdot 3^k$.  The marking $\phi$ of such a
minimal reduced marked graph 
is determined on $y_j, u_j$ and $v_j$ 
by a loop $w_j$ at $*$ and a single edge loop $\SP_j$ 
at $*$  in the graph as follows:
$$\begin{tabular}{lcl}
$\phi(y_j)$ & $=$ & $\bar w_j \star \SP_j \star  w_j$ \\
$\phi(u_j)$ & = & $\bar w_j \star e^j_1 \star \bar e^j_2 \star  w_j$\\
$\phi(v_j)$ & = & $\bar w_j \star e^j_2 \star \bar e^j_3 \star  w_j$, 
\end{tabular}$$
Note that the edge loops $\SP_j$ have to be different for $\phi$ to be a
homotopy equivalence.

By Theorem  $\ref{thm:initseg}$, both 
$\|\Lambda\|$ and $|\Lambda|$ are contractible.
So  Theorem \ref{thm:main} follows directly from the following statement:

\begin{thm}\label{thm:posetequivalence} There is a poset equivalence between
$L_{n,k}$  and $|\Lambda|$. 
\end{thm}

\begin{lem}\label{Lgkessential}
Let $e$ be an edge in a $(n,k)$-graph $(\Ga,\ul)$. Then there exists a
maximal forest in $\Ga-e$, that is a forest $\mathcal{F}$ such that 
$\Ga/\mathcal{F}$ is a (reduced) $(n,k)$-graph.
\end{lem}

\begin{proof}
Let $\hat \Ga$ be the graph of genus $n$ obtained from $\Ga$ by collapsing the cycles. Take a maximal forest $\hat\F$ in $\hat\Ga$. 
It can be extended to a maximal forest 
$\F$ of $\Ga$ by picking all but one edge in each cycle. 

 If $e$ is an edge in a cycle $\SP_j$, we
can choose such a maximal forest containing $\SP_j - \{e\}$ to obtain the
result.  If $e$ is not in a cycle, $e$ is non-separating in $\Ga$ by assumption and hence 
non-separating in $\hat\Ga$, so we can choose a forest $\hat\F$ in $\hat\Ga$ which does not contain $e$ 
and extend it to $\Ga$ as above. 
\end{proof}

\begin{proof}[Proof of the Theorem] 
Let $(\Ga,\phi)$ be a marked $(n,k)$-graph. 
Its underlying $(n,k)$-graph has $\SP_j=\phi(C_j)$ with
basepoint $\circ_j$,  and let $w_j$ be the path from $\circ_j$ to $*$ in $\Ga$,
image of the stem of the $j$th cycle in $G_{n,k}$, so that  
$\phi(y_j)=\bar w_j\star\SP_j\star w_j$.
Define a graph $\tilde\Ga$ 
by attaching $\theta^j$ to $\circ_j$ in $\Ga$ for each $j$ 
and define a marking $\tilde\phi:R_{n+3k}\to\tilde\Ga$ by 
$$\begin{tabular}{lcl}
$\tilde\phi(x_i)$ & $=$ & $\phi(x_i)$ \\
$\tilde\phi(y_j)$ & $=$ & $\bar w_j \star \SP_j \star  w_j$ \\
$\tilde\phi(u_j)$ & = & $\bar w_j \star e^j_1 \star \bar e^j_2 \star  w_j$\\
$\tilde\phi(v_j)$ & = & $\bar w_j \star e^j_2 \star \bar e^j_3 \star  w_j$. 
\end{tabular}$$
This defines an inclusion of $L_{n,k}$ into $\|\Lambda\|$.  Indeed, the
poset minimal elements of $L_{n,k}$ are mapped to $\Lambda$, so their stars
are in $\|\Lambda\|$. 

To show that it maps into the smaller
subcomplex $|\Lambda|$, we need to show that for each $(\Ga,\phi)$ in 
$L_{n,k}$, the pair $(\tilde\Ga,\tilde\phi)$ is essential in $\Lambda$, that is for each 
edge $e$ of $\tilde\Ga$, there is an invariant forest $\mathcal{F}\in\tilde\Ga-e$
such that $\tilde\Ga/\mathcal{F}\in\Lambda$. Note first that a maximal invariant forest of
$\tilde\Ga$ has to lie in $\Ga$. If $e$ lies in $\theta^j$ for some $j$,
choose a maximal forest $\mathcal{F}$ in $\Ga$ such that $\Ga/\mathcal{F}$ is
a reduced $(n,k)$-graph and hence $\tilde\Ga/\mathcal{F}$ is in $\Lambda$. 
On the other hand, if $e$ lies in $\Ga$, we know that such a forest exists in
$\Ga/e$ by Lemma \ref{Lgkessential}.

We claim that the image of $L_{n,k}$ equals $|\Lambda|$.
To show this, we need to show that the marking of any graph in $|\Lambda|$ is
induced by a marking from $G_{n,k}$ to the graph minus the theta subgraphs, 
that is we need to show how that the image of $y_j$ can be decomposed as
a cycle with a stem for each $j$ with at most a point in the intersection of any two 
different cycles,  
 and we need to show that the $u_j$'s and $v_j$'s are
mapped accordingly. 
This is true for the poset 
minimal elements, i.e. the elements of $\Lambda$: the
cycles are the $\SP_j$'s with basepoint $\circ_j=*$ and with stem  $w_j$.  
So we need to show that the cycle structure is carried by the admissible
blow-ups (that is the blow-ups which stay within $|\Lambda|$).

Let $(\Ga,\psi)\in |\Lambda|$ be a graph which has a cycle structure and let 
$\gamma$ be 
an ideal edge at a vertex $v$ of $\Ga$ such that $\Ga^\ga$ is still
in $|\Lambda|$.
We want to show that $\Ga^\ga$ has a cycle structure. 

\vs

\noindent
{\bf Case 1}\qua Suppose that $v\neq\circ_j$ for all $j$, so there is no
$\theta^j$ attached to $v$. Each cycle has at
most two edges ending at $v$. If $\gamma$ has either 0 or 2 half edges from $\SP_j$
for each $j$ except for at most one $j$, 
then the induced marking on $\Ga^\ga$ still has embedded
circles, with intersection possibly smaller than in $\Ga$. 

Suppose now that there is $j\neq j'$ such that both  $j$ and $j'$ have
exactly one half edge in $\ga$. Then the induced cycles in $\Ga^\ga$ are such
that $\SP_j\cap\SP_{j'}=\ga$, where $\ga$ is now the new edge in $\Ga^\ga$.
But there is no forest $\mathcal{F}$ in $\Ga^\ga-\{\ga\}$ such that
$\Ga^\ga/\mathcal{F}$ is in $\Lambda$. Indeed, if we do not collapse $\ga$
again, the induced marking will have the $j$th and $j'$th cycles intersecting
in $\ga$. Hence $\Ga^\ga$ is not in $|\Lambda|$.

\vs

\noindent
{\bf Case 2}\qua Suppose that there is at least one $j$ such that $v=\circ_j$. In particular,
there is at least one $\theta$-graph attached to $v$. Note first that the
equivariance conditions on $\ga$ imply that at most one half edge of a graph
$\theta^j$ may be in $\ga$, and if there is such an edge, the whole $\theta$-graph is
moved along the new edge $\ga$ in $\Ga^\ga$. 
Moreover, by the same argument as in Case 1, their can be at most one $j$ such
that $\SP_j$ has exactly one half edge in $\ga$. 

Suppose that there is a $j$ such that $\SP_j$ has two half edges in
$\ga$. Then $e_i^j$ must also be in $\ga$ for some $i$, otherwise 
any forest collapse which does not collapse $\ga$ 
will all induce a marking with
$\bar w_j\star\ga\star\SP_j\star\bar\ga\star w_j$ for $y_j$ but still 
$\bar w_j \star e^j_1 \star \bar e^j_2 \star  w_j$ for $u_j$. 
Similarly, if a half edge $e_i^j$ is in $\ga$ for some $j$, then the two half
edges of $\SP_j$ terminating at $v$ must also be in $\ga$. 
This shows that a blow up cannot disconnect $\theta_j$ from $\SP_j$. 

Lastly, if there is a $j$ such that only one half edge is in $\ga$, then the
induced marking will move $\circ_j$ and prolong $w_j$ if there is also an $e_i^j$
in $\Ga$. 
\end{proof}

\section{A space $L_n^k$ for the group $A_n^k$}\label{extra}

We recall from Section \ref{Anksec} that there is a central 
extension 
$$\Z^k\rar A_{n,k}\rar A_n^k,$$
where $A_n^k$ is a subgroup of $\Aut(F_{n+k})$.
To this group extension, corresponds  a sequence  of contractible spaces 
$$F_{n,k} \rar L_{n,k} \rar L_n^k$$
where the spaces are realizations of posets,  
the maps are poset maps and $F_{n,k}\cong \R^k$. 
We describe here the spaces $L_n^k$ and $F_{n,k}$.

\vs

As $A_n^k$ is a subgroup of $\Aut(F_{n+k})$, 
a natural way to construct a contractible space $L_n^k$ with an action of $A_n^k$ with finite stabilizers and finite quotient, is to consider a 
subcomplex of Auter space $X_{n+k}$, as small as possible to get as good an approximation as we can of the vcd. 
In the same spirit as the construction of $L_{n,k}$, one can consider the essential realization $|\Lambda'|$ of the initial segment $\Lambda'$ 
of reduced graphs of minimal norm with respect to the $n$-tuple of cyclic words $(y_1,\dots,y_k)$. 
We claim that $|\Lambda'|$ is homeomorphic to the following space $L_n^k$:

\begin{defn}
An {\em $\binom{k}{n}$-graph} $(\Gamma,\ul)$ is a graph $\Ga$ of genus $(n+k)$ with vertices of valence at least 3, except  
the basepoint which has valence at least 2, with no separating edges, and equipped with embeddings
 $l_j: S^1 \to \Gamma$ for $j=1,\ldots,k$ such that 
any two $\SP_j:=l_j(S^1)$ meet in at
most a point and the dual graph of the $\SP_j$'s is a forest.

A {\em marked $\nk$-graph} $(\Ga,\phi)$ is a graph $\Ga$ equipped with a pointed 
homotopy equivalence $\phi:G_{n,k} \to \Ga$ 
such that 
$(\Ga,\phi\circ \underline{\kappa})$ is an $\nk$-graph.

Two marked $\nk$-graphs $(\Gamma^1,\phi^1)$ and $(\Gamma^2,\phi^2)$ are
{\em equivalent} if there is a graph isomorphism 
$\psi: \Gamma^1 \to \Gamma^2$ 
such that $\psi\circ \phi^1(C_j)=\phi^2(C_j)$ 
and $\psi\circ \phi^1\simeq\phi^2$ (rel. $C_j$ for all $j$, up to a rotation). 
\end{defn}

Define now $L_n^k$ to be the poset of equivalence classes of 
marked $\nk$-graphs (or its realization) where the poset structure is given by the same admissible edge collapses 
as for $L_{n,k}$ (see p6).

\begin{thm}
There is a poset equivalence $L_n^k\to |\Lambda'|$ and hence $L_n^k$ is contractible.
\end{thm}

The proof is totally analogous to the proof of Theorem \ref{thm:posetequivalence}, although a little bit simpler as we do not need to 
carry around $\Theta$-graphs. The poset map $L_n^k\to X_{n+k}$ takes the minimal poset elements of $L_n^k$ to $\Lambda'$, and $L_n^k$ maps 
into the essential realization of $\Lambda'$ by Lemma \ref{Lgkessential}. On the other hand, we have $\Lambda'\subset L_n^k$ and the inclusion 
of $|\Lambda'|$ reduced to the case 1 in the proof of Theorem \ref{thm:posetequivalence}. Finally, $|\Lambda'|$ is contractible by 
Theorem \ref{thm:initseg} with $G$ the trivial group.

\begin{prop}
The quotient space $L_n^k/A_n^k:=QL_n^k$ has dimension $2n+2k-2$ if $n \geq 1$
and dimension $k-1$ if $n=0$.
\end{prop}

Again, the proof is essentially the same as for $QL_{n,k}$ (in Theorem~\ref{vcd}). 

Note that Collins constructed a space for $A_0^k$ which is isomorphic to our space $L_0^k$ 
(up to basepoint questions) \cite{c}. 
In this particular case 
when $n=0$, the dimension is equal to the vcd.

There is a forgetful map $L_{n,k}\to L_n^k$. Let $F_{n,k}$ be the fiber over the rose $(R_{n+k},id)$ with the $\Z^k$-action induced by 
the action of $A_{n,k}$ on $L_{n,k}$. 

\begin{prop}
There is a poset equivalence $F_{n,k}\to \Z^k\times(\Z/2)^k$, where the poset structure on $\Z^k\times(\Z/2)^k$ is defined by 
$(\bar z,\bar p)\le (\bar z',\bar p')$ iff for all $i$ we have 
\begin{tabular}{lll}
$p_i=p_i'$ & and & $z_i=z_i'$ or\\
$p_i<p_i'$ & and & $z_i=z_i'$ or $z_i=z_i'+1$.
\end{tabular}

Moreover, the realization of the poset $\Z^k\times(\Z/2)^k$ is homeomorphic to $\R^k$ and the $\Z^k$-action on $F_{n,k}$ translates 
to the canonical $\Z^k$-action on $\R^k$.
\end{prop}

\begin{proof}
Suppose $(\Ga,\phi)\in F_{n,k}$. Then $\Ga$ is the rose $R_{n+k}$ with additional basepoints for the cycles, which are either equal 
to the basepoint of the graph, or a valence two vertex on the petal. Up to $(n,k)$-graph isomorphisms, there are $(\Z/2)^k$ possibilities. We identify 
$\bar p\in (\Z/2)^k$ with the rose which has a valence two vertex on the $(n+j)$th petal if $p_j=1$ and no such vertex if $p_j=0$. 

Now $\phi:G_{n,k}\to R_{n+k}$ is homotopic to the identity if the cycles are allowed to rotate. So there are $\Z^k$ possibilities 
for $\phi$, where $\bar z\in \Z^k$ for a rose $(R,\bar p)$ represents the map $\phi$ which maps the cycles via the identity, according to the 
basepoints parametrized by $\bar p$, and wraps the $j$th stem $z_j$ times around the $j$th cycle. 

If $p_j=0$ for all $j$, then the graph is reduced. If there exists a $j$ such that $p_j=1$, then there are two possible collapses
on the $j$th petal, one which leaves 
$z_j$ constant, and one which increases $z_j$ by 1. (We choose the identification $\phi\rightsquigarrow \bar z$ so that this is the case.)
This corresponds to the poset structure of $\Z^k\times(\Z/2)^k$.

We are left to show that the realization of the poset $\Z^k\times(\Z/2)^k$ is homeomorphic to $\R^k$ and that the induced $\Z^k$ action is the canonical action.
Embed the set $\Z^k\times(\Z/2)^k$ in $\R^k$ by mapping $(\bar z,\bar p)$ to 
$(\bar z+\frac{\bar p}{2})=(z_1+\frac{p_1}{2},\dots,z_k+\frac{p_k}{2})\in \R^k$. 
Let $C_{\bar z}$ be the cube of length 1 in $\R^k$ with bottom corner $\bar z$ and top corner $\bar z+\bar 1$.
Let $c_{\bar z}=(z_1+\frac{1}{2},\dots,z_k+\frac{1}{2})$ be the center of the cube. Then we have  
$C_{\bar z}\cap F_{n,k}=\{(\bar z,\bar p)\in F_{n,k}| (\bar z,\bar p)\le c_{\bar z}=(\bar z,\bar 1)\}$. 
Similarly, the center of each face is the maximal element of the face in the poset. It follows that  
the realization of $\Z^k\times(\Z/2)^k$ is $\R^k$. The $\Z^k$-action, which is by addition on the $\bar z$-coordinate, is then by translation.
\end{proof}

The contractibility of $L_n^k$ and $F_{n,k}$ leads to an alternative proof of the contractibility of $L_{n,k}$. This second approach, 
suggested to us by the referee, makes more apparent the fact that the $\Theta$-graphs in the original proof of the contractibility 
of $L_{n,k}$ only play the role of adding basepoints.

\section{Presentation}\label{pressec}

We denote by $P_{i,j}$ the automorphism of the free group $F_n=\langle x_1,\dots,x_n\rangle$ which permutes $x_i$ and $x_j$ and leaves the other generators fixed, 
and we denote by $I_i$ the automorphism that maps $x_i$ to its inverse and fixes the other generators. 

For Whitehead moves, we use the following shortened notation: for $\e,\eta=\pm 1$,

$\bullet$\qua $(z_i^\e;z_j^\eta)$ denotes the automorphism which maps $z_i$ to $z_i z_j^\eta$ if $\e=1$ and to 
  $z_j^{-\eta}z_i$ if $\e=-1$, and fixes the other generators;

$\bullet$\qua  $(z_i^\pm;z_j^\eta)$denotes the automorphism which maps $z_i$ to $z_j^{-\eta}z_i z_j^\eta$ and fixes the other generators.

In terms of the discussion of Whitehead moves at the beginning of
Section \ref{fixedsection} before Theorem \ref{thm:initseg}, $(z_i^\e;z_j^\eta):=(\{z_i^\e,z_j^\eta\};z_j^\eta)$ and 
$(z_i^\pm;z_j^\eta):=(\{z_i,z_i^{-1},z_j^\eta\};z_j^\eta)$.

Note that $(z_i^\e;z_j^\eta)^{-1}=(z_i^\e;z_j^{-\eta})$ and  $(z_i^\pm;z_j^\eta)^{-1}=(z_i^\pm;z_j^{-\eta})$. 
We will use the second notation when writing relations for simplicity.

\begin{thm}\label{presAnk1}
The following gives a presentation for $\A_n^k$:

Generators: 

\begin{tabular}{ll}
$P_{i,j}$ & for  $1\le i,j \le n$ and $i\neq j$\\
$I_i$ &  for  $1\le i \le n$ \\
$(x_i^\e;z)$ & for  $1\le i \le n$, $\e=\pm 1$ and $x_i\neq z\in \{x_1,\dots,x_n,y_1,\dots,y_k\}$\\
$(y_i^\pm;z)$ &   for $1\le i \le k$ and $y_i\neq z\in \{x_1,\dots,x_n,y_1,\dots,y_k\}$\\
\end{tabular}

\noindent
Relations: For $z,z_i\in \{x_1,\dots,x_n,y_1,\dots,y_k\}$  and $w,w_i=x_{j_i}^{\delta_{i}}$ or $y_{j_i}^\pm$,

\begin{tabular}{ll}
Q1 & Relations in $\Aut(F_n)$ for $\{(x_i^{\e};x_j),P_{i,j},I_j\}$ \\
Q2 & $(w_1;z_1)(w_2;z_2)=(w_2;z_2)(w_1;z_1) \ \ \ $ for $w_1\neq w_2$ 
and $z_i^{\pm 1} \notin \{w_1,w_2\}$  \\
Q3 & $(1)\ \ (y_i^\pm;x_j)\ P_{j,l}=P_{j,l}\ (y_i^\pm;x_l)$ \\ 
   & $(2)\ \ (y_i^\pm;x_j)\ I_j=I_j\ (y_i^\pm;x_j^{-1})$  \\
Q4 & $(1)\ \ (w;x_j^\eta) (x_j^\eta;z) (w;x_j^{-\eta})=(w;z)(x_j^\eta;z)$ \\
   & $(2)\ \ (y_i^\pm;z^\e)(w;y_i)(y_i^\pm;z^{-\e})=(w;z^{-\e})(w;y_i)(w;z^\e)\ $ \\
Q5 & $(y_i^\pm;x_j^\eta)(x_j^{-\eta};y_i)=(x_j^\eta;y_i^{-1})(y_i^\pm;x_j^\eta)$
\end{tabular} 

whenever these symbols denote generators or their inverses.
\end{thm}

In particular, for $n=0$ we recover McCool's presentation  of the pure symmetric (or basis-conjugating) 
automorphisms of free groups \cite{mc4}, which has generators
$\{(y_i^\pm;y_j)\  |\ 1\le i,j\le k, i\neq j\}$ and relations

\begin{sam2}{\arabic}{Z}{}
\item $(y_i^\pm;y_j)(y_l^\pm;y_j)=(y_l^\pm;y_j)(y_i^\pm;y_j)$
\item $(y_i^\pm;y_j)(y_l^\pm;y_m)=(y_l^\pm;y_m)(y_i^\pm;y_j)$
\item $(y_i^\pm;y_j)(y_l^\pm;y_j)(y_i^\pm;y_l)=(y_i^\pm;y_l)(y_i^\pm;y_j)(y_l^\pm;y_j)$
\end{sam2}
where the indices are assumed to be different when given by different letters.

\begin{thm}\label{presAnk2}
The following gives a presentation for $\A_{n,k}$:

Generators: same as for $\A_n^k$ with $k$ additional generators $(y_i^\pm;y_i)$ for $1\le i\le k$.

Relations: Q1, Q2, Q3, Q4 (never allowing a symbol $(y_i^\pm;y_i)$), and 

\begin{tabular}{ll}
Q2\/$'$ & $(y_i^\pm;y_i)a=a(y_i^\pm;y_i)\ \ $ with $a=P_{i,j}, I_j$ or $(w;z)$\\
Q5\/$'$ & $(y_i^\pm;x_j^{\eta})(x_j^{-\eta};y_i)(y_i^\pm;x_j^{-\eta})(x_j^{\eta};y_i)=(y_i^\pm;y_i^{-1})$
\end{tabular}

whenever these symbols denote generators or their inverses.
\end{thm}

The generator $(y_i^\pm;y_i):=(x_1,\dots,x_n,1,\dots,y_i,\dots,1)\in \A_{n,k}$
 should be thought of as the conjugation of $y_i$ by itself, even though this is trivial as an automorphism 
of the free group.

\begin{proof}[Proof of Theorem \ref{presAnk1}]  
The group 
$\Aut(F_{n+k})$ acts on the set of cyclic words with letters in $L_{n+k}:=\{x_1,x_1^{-1},\dots,y_k,y_k^{-1}\}$ and 
 $\A_n^k$ is the stabilizer of the $k$-tuple of cyclic words $U=(y_1,\dots,y_k)$.  
Theorem 1.1 in \cite{mc3} says that such a group is finitely presented.
As remarked in  \cite[Section 4(1)]{mc3}, the proof actually  constructs a finite 2-dimensional complex $K$ whose fundamental 
group is this stabilizer and hence gives a presentation. 
We analyze this complex.

The vertices of the complex in our case are 
$$K_0:=\{V\ \textrm{minimal tuple }|V=UP\ \textrm{for some}\ P\in \Aut(F_{n+k})\}=\{U\s|\s\in\Omega_k\}$$
where a minimal tuple is a tuple of reduced cyclic words of minimal length in the orbit of $U$ and 
$\Omega_k\cong \Sigma_k\wr \Z_2$ is the subgroup of $\Aut(F_{n+k})$ of permutations and inverses of the $y_i$'s.
So $K_0\cong \Omega_k$.

Let  
 $\mathcal{W}_{n+k}$ denote the set of Whitehead automorphisms for $\Aut(F_{n+k})$.  We recall that $\mathcal{W}_{n+k}$ is the finite set 
$T_1\cup T_2$ where $T_1=\Omega_{n+k}$ is the set of permutations and inverses  and
 $T_2=\{(A;a)|A\subset L_{n+k},a\in A \ \textrm{and}\ a^{-1}\notin A\}$ with $(A;a)$ denoting the automorphism mapping $z\in L_{n+k}$
 to $za$ (and $z^{-1}$ to $a^{-1}z^{-1}$) if $a\neq z\in A$ and $z^{-1}\notin A$, 
mapping $z$ to $a^{-1}za$ (and $z^{-1}$ to $a^{-1}z^{-1}a$) if $z,z^{-1}\in A$ and fixing $z$ otherwise.
The set of edges in the 
complex is 
$$K_1:=\{(V_1,V_2;P)|P \in\mathcal{W}_{n+k} \ \textrm{and}\ V_1P=V_2\}/\sim.$$
Note that if $P\in \mathcal{W}$, then $P^{-1}$ is also in $\mathcal{W}$. 
The directed edge $e=(V_1,V_2;P)$, from $V_1$ to $V_2$, is identified with $\bar e=(V_2,V_1;P^{-1})$ in $K_1$. 
Finally the 2-cells are all the cells coming from the relations R1 to R10 in \cite{mc1}, attached to directed paths.
We recall these relations below.

There are two types of edges in $K$, coming from the two types of Whitehead automorphisms:

{\bf Type I edges}\qua $(U\s,U\s\tau;\gamma)$ for $\s,\tau\in \Omega_k$ and $\ga \in\Omega_{n+k}$ such that 
$(U\s)\gamma=U\s\tau$. In particular, $\gamma$ is only allowed to permute the $y_i$'s among themselves. Hence 
$\gamma\in \Om_n\times \Om_k\subset \Om_{n+k}$ and $\gamma=(\overline{\ga},\tau)$. 
So the set of type I edges between any two vertices is isomorphic to $\Om_n$.

{\bf Type II edges}\qua Type II edges do not involve permutations so they can exist only from one vertex to itself in our complex. 
Moreover, they have to preserve the cyclic words $y_1,\dots,y_k$, so they have the form 
$$(U\s,U\s;(X\sqcup\{y_{i_1},y_{i_1}^{-1},\dots,y_{i_r},y_{i_r}^{-1}\}\sqcup \{a\};a)$$
where $X\subset L_x=\{x_1,x_1^{-1},\dots,x_n,x_n^{-1}\}$  and $a\in L_{n+k}$.

Choose a maximal tree $T$ in $K_1$ with edges of type $(U\s,U\s\tau;\tau)$ for some $\tau\in\Om_k\cong\{1\}\times\Om_k$. 
The generators of $\pi_1{K}\cong \A_n^k$ are given by the edges in $K_1-T$. 

If we consider the subcomplex of $K$ formed by taking all the edges of the form  $(U\s,U\s\tau;\tau)$ for $\tau\in\Om_k$ and 
the 2-cells given by the relations R7 (which is a set of relations for $\Om_{n+k}$) restricted to these edges, we obtain a retraction of the 
2-skeleton of $E\Om_k$, which is simply connected. 
It follows that all edges of type  $(U\s,U\s\tau;\tau)$ give trivial generators in $\A_n^k$. 

Moreover, 
$ (U\s,U\s\tau;(\ga,\tau))\sim (U,U;(\ga,1))$
as $(1,\s)(\ga,\tau)=(\ga,1)(1,\s\tau)$ in $\Omega_{n+k}$ (giving thus a relation in R7) and
$ (U\s,U\s;(A;a))\sim (U,U;(A\s^{-1};a\s^{-1}))$
as $(1,\s)(A;a)=
(A\s^{-1};a\s^{-1})(1,\s)$ by R6. (R6 says that $T^{-1}(A;a)T=(AT;aT)$ for any $T\in T_1$.)
So the loops at the vertex $U$ generate the group. Also, the relations given by 2-cells in the complex all follow from the 
relations R1-R10 interpreted on the generators  $(U,U;\ga)$ and $(U,U;(A;a))$. Indeed, the other relations are conjugated to these 
(except for the relations in R6 and R7 already used to identify generators) 
and thus do not bring anything new by R6, used in the case $T=(1,\s)$. 

R1 identifies the inverses in $T_2$ and R2 says that $(A;a)(B;a)=(A\cup B;a)$ when $A\cap B=\{a\}$, showing thus that the generators we 
propose are indeed generators for the group.  Note that R2 also implies the commutation relations of Q2 when $z_1=z_2$.
As the relations in the Theorem are easily verified, we are left to show that they 
imply the relations R3-R10.

R3 is a set of commutation relations and follows from Q2 when $z_1\neq z_2$.

R4 says that  $(B;b)^{-1}(A;a)(B;b)=(A+B-b;a)$ when $A\cap B=\emptyset$, $a^{-1}\notin B$ and
$b^{-1}\in A$. We need to deduce it from our relations 
when $(A;a)$ and $(B;b)$ define loops at $U$ in $K$. In particular, we need $b=x_i^\e\in L_x$ as $b^{-1}\in A$ but $b\notin A$ and $b\neq a^{-1}$. 
Using Q2, we can isolate a minimal case on the left hand side:
 $(B-w;x_i^\e)^{-1}(w;x_i^{-\e})(x_i^{-\e};a)(w;x_i^\e)(A-x_i^{-\e};a)(B-w;x_i^\e)$.
We then use Q4 (1) and repeat the procedure with all the elements of $B$.

R5 says that $(A;a)(A-a+a^{-1};b)=\left( \begin{array}{cc} a & b \\ \downarrow & \downarrow \\ b^{-1} & a \end{array}\right) 
                     (A-b+b^{-1};a)$ 
when $b\in A$, $b^{-1}\notin A$ and $a\neq b$. In our case, we need to have $a,b\in L_x$. The relation follows from Q1 when only  
$x$'s are involved. Using R2, we isolate $(b;a)(a^{-1};b)$ in the left hand side and use R5 in this minimal case (that is use Q1). 
We then use 
Q3 to move the matrix (element of T1) in front of the left hand side and finally use R4.

R6, which was recalled above, follows from Q1 and Q3 (using Q2 as usual).

R7 gives a set of relations for $\Om_n$, which is included in Q1.

R8 says that $(A;a)=(L_{n+k}-a^{-1};a)(L_{n+k}-A;a^{-1})=(L_{n+k}-A;a^{-1})(L_{n+k}-a^{-1};a)$. This just follows 
from R1-R2.

R9 says that $(L_{n+k}-b;b^{-1})(A;a)(L_{n+k}-b^{-1};b)=(A;a)$ if $b,b^{-1}\notin A$.
In particular,  $a\neq b$. This relation follows from Q4. 
Indeed, if $a=x_j^\delta\in L_x$, isolate $(x_j^{-\delta};b^{-1})$ 
from of the left hand side and use Q4 (1) to pass it over $(A;x_j^\delta)$, which gives 
$(L_{n+k}-b-x_j^{-\delta};b^{-1})(A;x_j^\delta)(A-x_j^\delta+b^{-1};b^{-1})(x_j^{-\delta};b^{-1})(L_{n+k}-b^{-1};b)$.

Now simplify and isolate $(x_j^\delta;b)$ and use again Q4(1) to pass it over $(A;x_j^\delta)$, which gives
$(L_{n+k}-b-x_j^{-\delta};b^{-1})(x_j^{\delta};b)(A-x_j^\delta+b;b)(A;x_j^\delta)(L_{n+k}-A-x_j^{-\delta}-b^{-1};b)$.

This simplifies to $(A;x_j^\delta)$ as required.

If on the other hand $a=y_i^\e\in L_y$, isolating $(y_i^\pm;b^{-1})$ in the left hand side and using Q4 (2) to 
pass it over $(A;y_i^\e)$ gives  

$(L_{n+k}-b-y_i^\pm;b^{-1})(A-y_i^\e+b;b)(A;y_i^\e)(A-y_i^\e+b^{-1};b^{-1})(y_i^\pm;b^{-1})(L_{n+k}-b^{-1};b)$,

which simplifies to $(A;y_i^\e)$.

R10  says that 
$(L_{n+k}-b;b^{-1})(A;a)(L_{n+k}-b^{-1};b)=(L_{n+k}-A;a^{-1})$ if $b\neq a$, $b\in A$ and $b^{-1}\notin A$. So $b\in L_x$ in our case. 
If $a\in L_x$, we can follow \cite{mc1} and show that the relation follows from R9 and R5. 
We do the case  $a\in L_y$. So we take $a=y_j^\eta$, $b=x_i^\e$ 
The left hand side in the equation equals

 $(L_{n+k}-x_i^\e-y_j^\pm;x_i^{-\e})(y_j^\pm;x_i^{-\e})(x_i^\e;y_j^\eta)(A-x_i^\e;y_j^\eta)(L_{n+k}-x_i^{-\e};x_i^\e)$

which, using Q5, gives

 $(L_{n+k}-x_i^\e-y_j^\pm;x_i^{-\e})(x_i^{-\e};y_j^{-\eta})(y_j^\pm;x_i^{-\e})(A-x_i^\e;y_j^\eta)(L_{n+k}-x_i^{-\e};x_i^\e)$.

Now we use Q4 (1) to pass $(x_i^{-\e};y_j^{-\eta})$ to the left, and Q4(2) to pass $(y_j^\pm;x_i^{-\e})$ over $(A-x_i^\e;y_j^\eta)$, which gives

 $(x_i^{-\e};y_j^{-\eta})(L_{n+k}-x_i^\e-y_j^\eta;y_j^{-\eta})(L_{n+k}-x_i^\e-y_j^\pm;x_i^{-\e})(A-y_j^\eta;x_i^\e)(A-x_i^\e;y_j^\eta)$\nl
$(A-x_i^\e-y_j^\eta+x_i^{-\e};x_i^{-\e}) (y_j^\pm;x_i^{-\e})(L_{n+k}-x_i^{-\e};x_i^\e)$.

This simplifies to 
$(L_{n+k}-A;y_j^{-\eta})$ as required.
\end{proof}

To prove Theorem \ref{presAnk2}, we need a lemma about group extensions. 

Let $1\rar A\stackrel{i}{\rar} E\stackrel{\pi}{\rar} G\to 1$ be a central extension of a group $G$ by an abelian group $A$, with a 
 section of sets $s:G\to E$ (so $\pi\circ s=id$) satisfying $s(1)=1$ and $f:G\times G\to A$ the classifying cocycle of the extension, defined by 
$$s(g)\ s(h)=i\!\circ\! f(g,h)\ s(gh).$$
Note that if $r_1\dots r_m=1$ in $G$, then $s(r_1\dots r_m)=1$ in $E$
and hence
$$s(r_1)\dots s(r_m)=f(r_1,r_2)f(r_1r_2,r_3)\dots f(r_1\dots r_{m-1},r_m).$$ 
We denote by 
$s^f$ the operator
lifting the relations of $G$ to relations in $E$ this way.

\begin{lem}
Let $1\rar A\stackrel{i}{\rar} E\stackrel{\pi}{\rar} G\to 1$ be a central extension as above.
If $A$ and $G$ are finitely presented, with presentation 
\begin{gather*}A=\langle a_1,\dots,a_k\rangle/R_A\\ 
G=\langle g_1,\dots,g_n\rangle/R_G.\end{gather*} 
Then $E$ is finitely presented, with presentation 
$$E=\langle i(a_1)\dots,i(a_k),s(g_1),\dots,s(g_n)\rangle/[A,G]\cup i(R_A)\cup s^f(R_G)$$
where $[A,G]$ is the set of commutation relations $i(a_j)s(g_i)=s(g_i)i(a_j)$ and 
$s^f(R_G)$ is the lift of the relations $R_G$ to $E$ twisted by $f$ as defined above. 
\end{lem}

\begin{proof}
Consider the free group $F:=\langle i(a_1)\dots,i(a_k),s(g_1),\dots,s(g_n)\rangle$. The natural map $h:F\to E$ is surjective as $i\times s: A\times G\to E$ 
is a bijection of sets (since $s(G)$ is a set of coset representatives). The relations $[A,G]$, $i(R_A)$ and $s^f(R_G)$ are clearly 
satisfied in $E$, and so define elements in the kernel of $h$. 

Let $e=i(a_{j_1})\dots i(a_{j_p})s(g_{i_1})\dots s(g_{i_q})\in F$ and suppose that $h(e)=1$. 
(It is enough to consider this case as $[A,G]$ is satisfied in $E$.)
Then $\pi\circ h(e)=1$ in $G$, that is  
$g_{i_1}\dots g_{i_q}=1$. This relation in $G$ lifts to a  relation  $s(g_{i_1})\dots s(g_{i_q})=i(a)$ in $E$ for some $a\in A$. 
So $i(a_{j_1})\dots i(a_{j_p})s(g_{i_1})\dots s(g_{i_q})=i(a_{j_1})\dots i(a_{j_p})i(a)$. As $h(e)=1$, we have 
$a_{j_1}\dots a_{j_p}a=1$ is a relation in $A$ (as $i$ is an injective homomorphism). So the  
relation $h(e)=1$ follows from a relation in $s^f(R_G)$ and a relation in $i(R_A)$.
\end{proof}

\begin{proof}[Proof of Theorem \ref{presAnk2}]  
There is a central extension 
$1\to \Z^k \to \A_{n,k} \to \A_n^k \to 1$. By the Lemma, the generators of $\A_{n,k}$ are the union of the generators of $\A_n^k$ 
with those of $\Z^k$. The last generators are denoted $(y_i^\pm;y_i)$, for $1\le i\le k$, in the Theorem. Then relations for 
$\A_{n,k}$ are given by the relations for $\Z^k$, the commutations $[\Z^k,\A_n^k]$ and the lifts of the relations from $\A_n^k$. 
Q2' takes care of the first two sets of relations and  Q1,Q2,Q3,Q4 and Q5' is the lift of the relations from $\A_n^k$.
\end{proof}

In \cite{bl}, the group $\A_0^k$ 
is studied (denoted $C_k$ in the paper). 
They show that 
$$H_1(\A_0^k)=\bigoplus_{1\le i,j\le k,\ i\neq j}\Z[y_i^\pm;y_j] $$
where $\Z[y_i^\pm;y_j]$ is the free abelian group generated by $[y_i^\pm;y_j]$. 
Recall that $\A_0^k$ can be thought of as an extended pure braid group (which is the mapping class group of a punctured sphere). 
Note in particular that, as for the pure braid group, $H_1$ 
is not stable with respect to $k$. 
The presentation given above allows us to calculate $H_1$ in the higher genus cases. 

\begin{thm}\label{H1}
$H_1(\A_{n,k})\cong \Z_2 \cong H_1(\A_n^k)$ if $n>2$
\end{thm}

\begin{proof}
If $n>2$, we have $H_1(\Aut(F_n))=\Z_2$ generated by $[P_{i,j}]=[I_k]$ (the generators $(x_i^\delta;x_j)$ being trivial in the abelianization). 

For $\A_n^k$, Q4 (1) implies that $(w;z)$ is trivial in the abelianization for any $w=x_i^\delta$ or $y_i^\pm$ and $z\in\{x_1,x_2,\dots,y_k\}$.
 The only relations involving $P_{i,j}$ or $I_k$ which are not in $\Aut(F_n)$ give trivial relations in 
the abelianization, hence $H_1(\A_n^k)=\Z_2$. 

For $\A_{n,k}$, the additional relation Q5' implies that the generators $(y_i^\pm;y_i)$ are also trivial in the abelianization and hence we get the same result. 
\end{proof}

\Addresses\recd


\begin{thebibliography}{1}

\bibitem{bmmm}
N. Brady, J. McCammond, J. Meier, and A. Miller,
{\em The pure symmetric automorphisms of a free group form a duality group},
J. Algebra 246 (2001) 881-896.

\bibitem{b}
K. Brown,
{\em Cohomology of Groups},
Springer-Verlag, New York-Berlin, 1982.  

\bibitem{bl}
A. Brownstein and R. Lee, {\em Cohomology of the group of motions of 
$n$ strings in 3-space}, Contemp. Math. vol. 150 (1993) 51-61. 


\bibitem{c}
D. J. Collins,
{\em Cohomological dimension and symmetric automorphisms of a
free group},
Comment. Math. Helv. 64 (1989) 44-61.

\bibitem{cv}
M. Culler and K. Vogtmann,
{\em Moduli spaces of graphs and automorphisms of free groups},
Invent. Math. 84 (1986) 91-119.

\bibitem{d} D. Dahm, {\em A generalization of braid theory}, Ph.D. Thesis, Princeton Univ., 1962.

\bibitem{g} D. Goldsmith, {\em The theory of motion groups},  Mich. Math. J. 28 (1981), 3-17.

\bibitem{h} J. L. Harer, {\em Stability of the homology of the mapping 
class groups of orientable surfaces}, Annals Math. 121 (1985), 
215-249.

\bibitem{H} A. Hatcher, {\em Algebraic Topology}, Cambridge University Press, 2002.

\bibitem{hv}
A. Hatcher and K. Vogtmann,
{\em Cerf theory for graphs},
J. London Math. Soc. (2) 58 (1998) 633-655.

\bibitem{hw} A. Hatcher and N. Wahl, {\em Stabilization for the automorphisms of free groups with boundary}, preprint, \arxiv{math.GT/0406277}.

\bibitem{j0} 
C. A. Jensen,
{\em Cohomology of $\Aut(F_n)$ in the $p$-rank two case}, 
J. Pure Appl. Algebra 158 (2001) 41-81.

\bibitem{j1}
C. A. Jensen, {\em Contractibility of fixed point sets of auter space},
Topology Appl. 119 (2002) 287-304.

\bibitem{jmm}
C. A. Jensen, J. McCammond and J. Meier,
{\em The Euler characteristic of the Whitehead automorphism group of a free product},
preprint.

\bibitem{j2}
C. A. Jensen,
{\em Homology of holomorphs of free groups},
J. Algebra 271 (2004) 281-294.

\bibitem{k}
S. Krsti\'c,
{\em Actions of finite groups on graphs and related automorphisms
of free groups},
J. Algebra 124 (1989) 119-138.

\bibitem{kv}
S. Krsti\'c and K. Vogtmann,
{\em Equivariant outer space and automorphisms of
free-by-finite groups},
Comment. Math. Helv. 68 (1993) 216-262.

\bibitem{l} G. Levitt,  {\em Automorphisms of hyperbolic groups and graphs of groups}, preprint, \arxiv{math.GR/0212088}. 

\bibitem{mt} I. Madsen and U. Tillmann, {\em The stable mapping class group 
and $\mathcal{Q}(\mathbb{C}P^\infty_+)$}, Invent. math. 145 (2001), 509-544. 
 


\bibitem{mw} I. Madsen and M. Weiss, {\em The stable moduli space of Riemann surfaces: Mumford's conjecture}, preprint, \arxiv{math.AT/0212321}.

\bibitem{m}
W. Magnus, A. Karrass and D. Solitar, {\em Combinatorial group theory. 
Presentations of groups in terms of generators and relations}, 
Second edition, Dover Pub., New York, 1976.

\bibitem{mm2}
J. McCammond and J. Meier,
{\em The hypertree poset and the $L^2$ Betti numbers of the motion group of the trivial link},
Math. Ann. 328 (2004) 633-652.

\bibitem{mc1}
J. McCool, {\em A presentation for the automorphism group of a free group of finite rank}, J. London Math. Soc. (2) 8 (1974),259-266.

\bibitem{mc3}
J. McCool, {\em Some finitely presented subgroups of the automorphism group of a free group}, J. Algebra 35 (1975), 205-213.

\bibitem{mc4}
J. McCool, {\em On basis-conjugating automorphisms of the free groups}, Can. J. Math., 38 (1986) no. 6, 1525-1529.


\bibitem{mm}
D. McCullough and A. Miller,
{\em Symmetric Automorphisms of Free Products},
Mem. Amer. Math. Soc. 122 (1996), no. 582.



\bibitem{q}
D. Quillen,
{\em Homotopy properties of the poset of non-trivial $p$-subgroups of a 
group},
Advances in Math. 28 (1978) 101-128.


\bibitem{t} U. Tillmann, {\em On the homotopy of the stable mapping 
class group}, Invent. Math. 130 (1997), 257-275.


\bibitem{w}
N. Wahl,
{\em From mapping class groups to automorphisms of free groups}, preprint,
\arxiv{math.AT/0406278}.

\end{thebibliography}
\end{document}